\documentclass[11pt]{amsart}
\usepackage{amssymb}
\usepackage{amsbsy, amsthm, amsmath,amstext, amsopn}
\usepackage[all]{xy}
\usepackage{amsfonts}
\usepackage{amscd}
\usepackage{stmaryrd}
\usepackage{epigraph}
\usepackage{hyperref}
\usepackage{verbatim}
\usepackage{paralist}
\usepackage{yhmath}
\usepackage{graphicx}
\usepackage{sidecap}
\usepackage{epstopdf}
\newtheorem{thm}{Theorem} [section]

\newtheorem{lemma}[thm]{Lemma}

\newtheorem{corollary}[thm]{Corollary}
\newtheorem{prop}[thm]{Proposition}

\newtheorem{fundlemma}{Fundamental Lemma} [section]

\theoremstyle{definition}

\newtheorem{defn}[thm]{Definition}

\newtheorem{example}[thm]{Example}

\newtheorem{remark}[thm]{Remark}

\theoremstyle{remark}

\numberwithin{equation}{section}

\newcommand{\nc}{\newcommand}
\nc{\BA}{\mathbb A}
\nc{\BF}{\mathbb F}
\nc{\BH}{\mathbb H}
\nc{\BQ}{\mathbb Q}
\nc{\BN}{\mathbb N}
\nc{\BP}{\mathbb P}
\nc{\BR}{\mathbb R}
\nc{\BC}{\mathbb C}
\nc{\BZ}{\mathbb Z}

\nc{\CA}{\mathcal A}
\nc{\CB}{\mathcal B}
\nc{\CE}{\mathcal E}
\nc{\CF}{\mathcal F}
\nc{\CG}{\mathcal G}
\nc{\CK}{\mathcal K}
\nc{\CL}{\mathcal L}
\nc{\CM}{\mathcal M}
\nc{\CN}{\mathcal N}
\nc{\CO}{\mathcal O}
\nc{\CP}{\mathcal P}
\nc{\CS}{\mathcal S}
\nc{\CV}{\mathcal V}
\nc{\CX}{\mathcal X}

\nc{\fb}{\mathfrak b}
\nc{\fg}{\mathfrak g}
\nc{\fh}{\mathfrak h}
\nc{\fk}{\mathfrak k}
\nc{\fl}{\mathfrak l}
\nc{\fo}{\mathfrak o}
\nc{\fp}{\mathfrak p}
\nc{\fs}{\mathfrak s}
\nc{\ft}{\mathfrak t}

\nc{\on}{\operatorname}
\nc{\ol}{\overline}
\nc{\wh}{\widehat}

\nc{\Hom}{\on{Hom}}
\nc{\End}{\on{End}}
\nc{\Aut}{\on{Aut}}
\nc{ \Tr}{\on{Trace}}
\nc{\Id}{\on{Id}}
\nc{\Ind}{\on{Ind}}
\nc{\Vect}{\on{Vect}}
\nc{\Spec}{\on{Spec}}
\nc{\Bun}{\on{Bun}}
\nc{\Out}{\on{Out}}
\nc{\Inn}{\on{Inn}}

\nc{\bs}{\backslash}

\nc{\adquot} {/}
\begin{document}

\title{The Geometric Nature of the Fundamental Lemma}

\author{David Nadler}
\address{Department of Mathematics\\Northwestern University\\Evanston, IL 60208-2370}
\email{nadler@math.northwestern.edu}

\maketitle

\begin{abstract}
The Fundamental Lemma is a somewhat obscure combinatorial identity
introduced by Robert P. Langlands~\cite{tracesstable} as an ingredient in the theory of automorphic representations.
After many years of deep contributions by mathematicians working in representation theory, number theory, algebraic geometry,
and algebraic topology, a proof of the Fundamental Lemma was recently completed by Ng\^{o} Bao Ch\^au~\cite{ngoFL},
for which he was awarded a Fields Medal.
Our aim here is to touch on some of the beautiful ideas contributing to the Fundamental Lemma and its proof.
We highlight the geometric nature of the problem which allows one to attack a question in $p$-adic analysis with the tools of algebraic geometry. 
\end{abstract}

%

\setcounter{tocdepth}{1}
\tableofcontents


\section{Introduction}

Introduced by Robert P. Langlands in his lectures~\cite{tracesstable},
the Fundamental Lemma 
is a combinatorial identity which just as well could have achieved no notoriety. Here is Langlands commenting
on~\cite{tracesstable} on the IAS website~\cite{langlandsendoscopypage}:

\begin{quote}
``...the fundamental lemma which is introduced in these notes, is a precise and purely combinatorial statement that I thought must therefore of necessity yield to a straightforward analysis. This has turned out differently than I foresaw."
\end{quote}

Instead, the Fundamental Lemma has taken on a life of its own. Its original scope involves distributions on groups over local fields
($p$-adic and real Lie groups). Such distributions naturally arise as the characters of representations,
and are more than worthy of study on their own merit.
But with the immense impact of Langlands' theory of automorphic and Galois representations, many potential advances turn out to be downstream of the Fundamental Lemma. 
In particular, in the absence of proof, it became ``the bottleneck limiting progress on a host of arithmetic questions."~\cite{harrisintro}
It is rare that popular culture recognizes the significance of a mathematical result, much less an esoteric lemma, but the 
recent proof of the Fundamental Lemma by Ng\^o Bao Ch\^au~\cite{ngoFL}, for which he was awarded a Fields Medal, ranked seventh on {\it Time} magazine's
{\em Top 10 Scientific Discoveries of 2009} list.\footnote{
Teleportation was eighth.}

\medskip

\begin{figure}[h!] \label{orbits}

\includegraphics[scale=0.5]{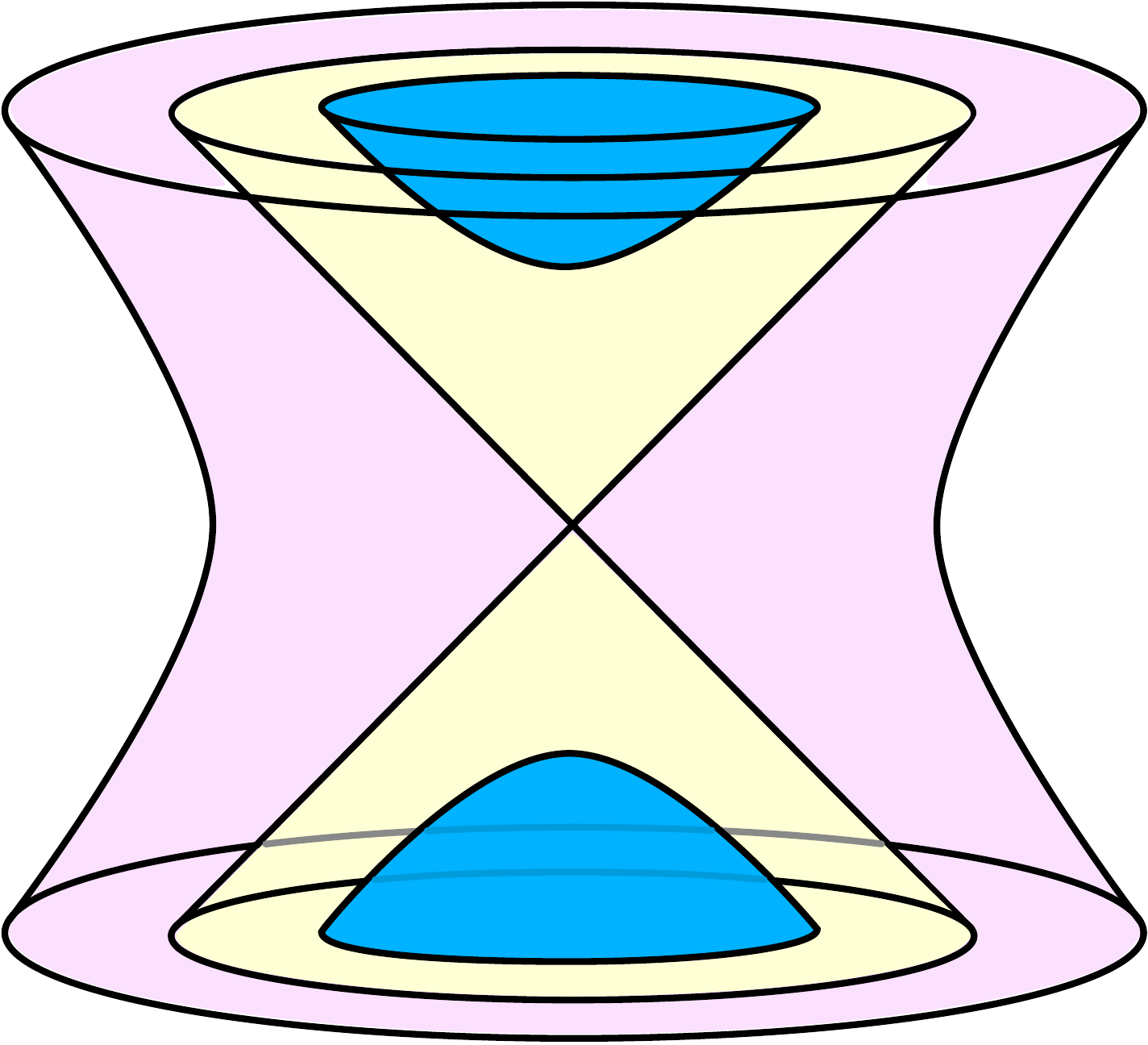}

\caption{Orbits of $SL(2, \BR)$ acting on its Lie algebra $\fs\fl(2, \BR) \simeq \BR^3$.}
 
\end{figure}

Before continuing, it might be useful to have in mind a  cartoon of the problem which the Fundamental Lemma solves. 
(In fact, what we present is an example of the case of real Lie groups resolved long ago by D. Shelstad~\cite{shelstad}.) 
Figure~\ref{orbits} depicts representative orbits for 
the real Lie group  $SL(2, \BR)$  
acting by conjugation on its Lie algebra $\fs\fl(2, \BR) \simeq \BR^3$ of traceless $2\times 2$ real matrices $A$.

Reading Figure~\ref{orbits} from outside to inside, one encounters
three types of orbits (hyperbolic, nilpotent, and elliptic) classified by the respective values of the determinant ($\det(A) < 0$, $\det(A) = 0$, and $\det(A)> 0$).
We will focus on the
two elliptic orbits $\CO_{A_+}, \CO_{A_-} \subset \fs\fl(2, \BR)$ through the 
elements
$$
A_+ = \left[
\begin{array}{cc}
0 & 1\\
 -1 & 0
\end{array}
\right],
\qquad
A_- = \left[
\begin{array}{cc}
0 &  -1\\
 1 & 0
\end{array}
\right].
$$
For a smooth compactly supported function $\varphi:\fs\fl(2, \BR) \to \BC$,  consider the distributions 
given by integrating
over the elliptic orbits
$$
\CO_{A_+}(\varphi) = \int_{\CO_+} \varphi,
\qquad
\CO_{A_-}(\varphi) = \int_{\CO_-} \varphi
$$
with respect to an invariant measure.

Observe that the two-dimensional complex vector space spanned by these distributions
admits the alternative basis
$$
\CO_{st}= \CO_{A_+} + \CO_{A_-},
\qquad
\CO_{tw} =   \CO_{A_+} - \CO_{A_-}.
$$

The first $\CO_{st}$ is nothing more than the integral over the union
$\CO_{A_+} \sqcup \CO_{A_-},
$ 
which  is the algebraic variety given by the equation $\det(A) =1$. 
%
It is called a {\em stable distribution} since the equation $\det(A)=1$ makes no reference to the field of real numbers $\BR$. 
Over the algebraic closure of complex numbers $\BC$, the equation $\det(A)=1$ cuts out  a single conjugacy class. In particular, $A_+$ and $A_-$ are both conjugate to the matrix
$$
 \left[
\begin{array}{cc}
i & 0\\
 0 & -i
\end{array}
\right]
\in \fs\fl(2, \BC).
$$
Thus the stable distribution can be thought of as an object of algebraic geometry rather than harmonic analysis on a real Lie algebra.

Unfortunately, there is no obvious geometric interpretation for the second $\CO_{tw}$.
(And one might wonder whether such a geometric interpretation could exist: the symmetry of  switching the terms of $\CO_{tw}$ gives
its negative.) It is called a {\em twisted distribution} since it is a sum of $\CO_{A_+}$ and $\CO_{A_-}$ with nonconstant coefficients. By its very definition, $\CO_{tw}$ distinguishes between the orbits  $\CO_{A_+}$ and $\CO_{A_-}$ though
there is not an invariant polynomial function which separates them. 
Indeed, as discussed above, over the complex numbers, 
they coalesce into a single orbit.

Langlands's theory of endoscopy, and in particular, the Fundamental Lemma  at its heart, confirms that indeed one can systematically express such twisted distributions in terms
of stable distributions. A hint more precisely: to any twisted distribution, there is assigned a stable distribution,
and to any test function, a transferred test function, such that the twisted distribution evaluated on the
original test function
is equal to the stable distribution evaluated on the transferred test function.

Detailed conjectures organizing the intricacies of the transfer of test functions first appear in Langlands's joint work with 
D. Shelstad~\cite{langlandsshelstad}. The shape of the conjectures for $p$-adic groups were 
decisively influenced by what could be more readily understood for real groups (ultimately building on work of 
Harish Chandra).
As Langlands and Shelstad note,``if it were not that [transfer factors] had been proved to exist over the real field~\cite{shelstad}, it would have been difficult to maintain confidence in the possibility of transfer or in the usefulness of endoscopy."

The  extraordinary difficulty of the Fundamental Lemma, and also its mystical power,  emanates from the fact that the sought-after  stable distributions live on the Lie algebras of groups
with little apparent relation to the original group. 
Applied to the example at hand, the general theory relates the  twisted distribution $\CO_{tw}$  to a stable distribution
on the Lie algebra $\fs \fo(2, \BR) \simeq \BR$ of the rotation subgroup $SO(2, \BR) \subset SL(2, \BR)$ which stabilizes $A_+$
or equivalently $A_-$. Outside of bookkeeping, this is empty of content since $SO(2, \BR)$ is abelian, and so its orbits in $\fs\fo(2, \BR)$ are single points. But the  general theory is deep and elaborate and leads to  surprising  identities of which the Fundamental Lemma
is the most basic and important.

\medskip


It should be pointed out that in the absence of a general proof, many special cases of the Fundamental Lemma
were established to spectacular effect. 
To name the most prominent applications without attempting to explain any of the terms involved, 
the proof of Fermat's Last Theorem depends upon base change for $GL(2)$, and ultimately
the Fundamental Lemma for cyclic base change~\cite{basechange}.
The proof of the local Langlands conjecture for $GL(n)$, a parameterization of the representations
of the group of invertible matrices with $p$-adic entries, depends upon automorphic induction, and ultimately the Fundamental Lemma for $SL(n)$ established by Waldspurger~\cite{waldsln}.

If the Fundamental Lemma had admitted an easy proof, it would have merited little mention in a discussion of these results.
But for the general theory of automorphic representations and Shimura varieties, ``...its absence rendered progress almost impossible for more than twenty years."~\cite{langlandslettertolang}
Its recent proof has opened the door to many further advances.
Arthur~\cite{arthurlocal, arthurreport} has outlined a program to obtain the local Langlands correspondence
for quasi-split classical groups from that of $GL(n)$ via twisted endoscopy
(generalizing work of Rogawski~\cite{rog} on the unitary group $U(3)$). 
In particular, this will provide a parameterization of the representations of orthogonal and symplectic matrix groups with $p$-adic entries.
Shin~\cite{shin} has constructed Galois representations corresponding to automorphic representations
in the cohomology of compact Shimura varieties, establishing the Ramanujan-Petersson conjecture 
for such representations. This completes earlier work of
Harris--Taylor,  Kottwitz and Clozel, and the Fundamental Lemma plays a key role in this advance. As Shin notes,
``One of the most conspicuous obstacles was the fundamental lemma, which had only been known in some special cases. Thanks to the recent work of Laumon-Ng\^o (\cite{LN}), Waldspurger (\cite{waltransfer}, \cite{walchar}, \cite{waltordu}) and 
Ng\^o (\cite{ngoFL}) the fundamental lemma (and the transfer conjecture of Langlands and Shelstad) are now fully established. This opened up a possibility for our work."
Morel~\cite{morel} has obtained similar results via a comprehensive study 
of the cohomology of noncompact Shimura varieties.

%
%
%
%

Independently of its applications, the peculiar challenge of the Fundamental Lemma has spurred many ingenious advances
of  standalone interest.
Its recent proof, completed by Ng\^o, spans many areas, appealing to remarkable new ideas in representation theory, model theory, 
algebraic geometry, and algebraic topology. 
A striking aspect of the story is its diverse settings. 
The motivation for and proof of the Fundamental Lemma sequentially cover the following algebraic terrain:
$$
\boxed{
\begin{minipage}{12.5cm}
\begin{center}
number fields
$\to$
$p$-adic fields
$\to$
Laurent series fields
$\to$
function fields
\end{center}
\end{minipage}
}
$$

One begins with a number field and the arithmetic problem of
comparing the anisotropic part of the Arthur-Selberg trace formula for different groups. This leads to the combinatorial question of the Fundamental Lemma about integrals over $p$-adic groups. Now not only do the integrals make sense for any local field,
but it turns out that they are  independent of the specific local field, and in particular its characteristic.
Thus one can work in the geometric setting of Laurent series with integrals over loop groups (or synonymously,
affine Kac-Moody groups). 
Finally, one returns to a global setting and performs analogous integrals along the Hitchin fibration for groups over the function field of a projective curve. In fact, one can interpret the ultimate results as precise geometric analogues
of the stabilization of the original trace formula. To summarize, within the above algebraic terrain, the main quantities to be calculated and compared for different groups are the following:
$$
\boxed{
\begin{minipage}{14cm}
\begin{center}
global anisotropic orbital integrals
$\to$
local twisted orbital integrals
$\to$\\
cohomology of affine Springer fibers $\to$
cohomology of anisotropic Hitchin fibers
\end{center}
\end{minipage}
}
$$

Over the last several decades, geometry has accrued a heavy debt to harmonic analysis and number theory: much of the representation
theory of Lie groups and quantum algebras, as well as gauge theory on Riemann surfaces and complex surfaces,
is best viewed through a collection of analogies (called the Geometric Langlands program,
and pioneered by Beilinson and Drinfeld) with Langlands's theory of automorphic and Galois representations.
Now with Ng\^o's proof of the Fundamental Lemma, and its essential use of loop groups and the Hitchin fibration, geometry has 
finally paid back some of this debt.

In retrospect, one is led to the question: {\em Why does geometry play a role in the Fundamental Lemma?} Part of the answer is implicit
in the fact that the $p$-adic integrals involved turn out to be characteristic independent. 
They are truly of motivic origin, reflecting universal polynomial equations rather than
analysis on the $p$-adic points of groups. Naively, one could think of the comparison between counting matrices
over a finite field with a prescribed characteristic polynomial versus counting those with prescribed eigenvalues.
More substantively, one could keep in mind Lusztig's remarkable construction of the characters of finite simple groups of Lie type
(which encompass ``almost all" finite simple groups).
Though such groups have no given geometric structure, their characters can be
uniformly constructed by recognizing the groups are the solutions to algebraic equations.
For example, though we may care primarily about characters of $SL(n, \BF_p)$, it is important to think not only about the set of such matrices, but also the determinant equation  $\det(A)= 1$
which cuts them out.

In the case of the Fundamental Lemma, the Weil conjectures ultimately imply that the $p$-adic integrals to be evaluated are shadows of the cohomology of 
algebraic varieties, specifically the affine Springer fibers of Kazhdan-Lusztig. 
Therefore one could hope to apply the powerful tools of mid-to-late 20th century algebraic geometry -- such as Hodge theory, Lefschetz techniques, sheaf theory, and homological algebra -- in the
tradition pioneered by Weil, Serre, Grothendieck, and Deligne.
One of Ng\^o's crucial, and possibly indispensable, contributions is to recognize that the technical structure needed to proceed, in particular the purity of the Decomposition Theorem of Beilinson-Bernstein-Deligne-Gabber, could be found in a return to the global setting
of the Hitchin fibration.

\medskip

Our aim in what follows is to sketch some of the beautiful ideas contributing to the Fundamental Lemma and its proof. 
The target audience is not experts in any of the subjects discussed but rather mathematicians interested in having some sense of the essential seeds from which a deep and intricate theory flowers.
We hope that in a subject with great cross-pollination of ideas, 
an often metaphorical account of important structures could prove useful.

Here is an outline.

In Section~\ref{secttrace} immediately below, 
we recall some basics about characters of representations and number fields leading to 
a very rough account of the Arthur-Selberg  trace formula for compact quotient.

In Section~\ref{sectstable}, 
we introduce the problem of stability of conjugacy classes, the twisted orbital integrals and endoscopic groups
which arise, and finally arrive at a statement of the Fundamental Lemma.

In the remaining sections, we highlight some of the beautiful mathematics of broad appeal which either contribute to the proof of the Fundamental
Lemma or were spurred by its study.
Some were invented specifically to attack the Fundamental Lemma, while others have their own pre-history but
are now inextricably linked to it.\footnote{Like Tang to NASA.}

In Section~\ref{sectlocal}, 
we introduce the affine Springer fibers and their cohomology which are the motivic avatars of the orbital integrals of the Fundamental Lemma.
We then discuss the equivariant localization theory of Goresky-Kottwitz-MacPherson developed to attack 
the Fundamental Lemma. Strictly speaking, it is not needed for Ng\^o's ultimate proof, but it both set the scene
for much of Laumon and Ng\^o's further successes, and has inspired an entire industry of combinatorial geometry.

In Section~\ref{sectglobal}, 
we summarize and interpret several key aspects of Ng\^o's proof of the Fundamental Lemma.
In particular, we discuss Laumon's initial forays to a global setting, and Ng\^o's Support Theorem which
ultimately provides the main technical input.

Finally, in Section~\ref{sectfuture}, we discuss some directions for further study.

\medskip

There are many precise and readable, though of necessity lengthy,  accounts of the mathematics we will discuss. In particular, 
we recommend the reader read everything by Langlands, Kottwitz, and Arthur, and time permitting, read all of 
Drinfeld and Laumon's lecture notes such as~\cite{drinfeldnotes, laumoncmi, laumonnewton}. For the Fundamental Lemma and its immediate neighborhood, there are Ng\^o's original paper~\cite{ngoFL}, and the long list of references therein,
Hales's beautifully concise paper~\cite{halessummary}, and the immensely useful book project organized by Harris~\cite{harrisintro}.

%
%
%
%
%
%

%
%
%
%
%

%

%
%

\subsection{Acknowledgements}
This document is intended as a writeup of my anticipated  talk at the Current Events Bulletin Session
 at the Joint Mathematics Meetings in New Orleans, January 2011.
I would like to thank the committee and organizers for the impetus
to prepare this document.

I am particularly indebted to D. Ben-Zvi for his generosity in sharing ideas during many far-ranging discussions.

\begin{figure}[h!] \label{sp4fig}

\includegraphics{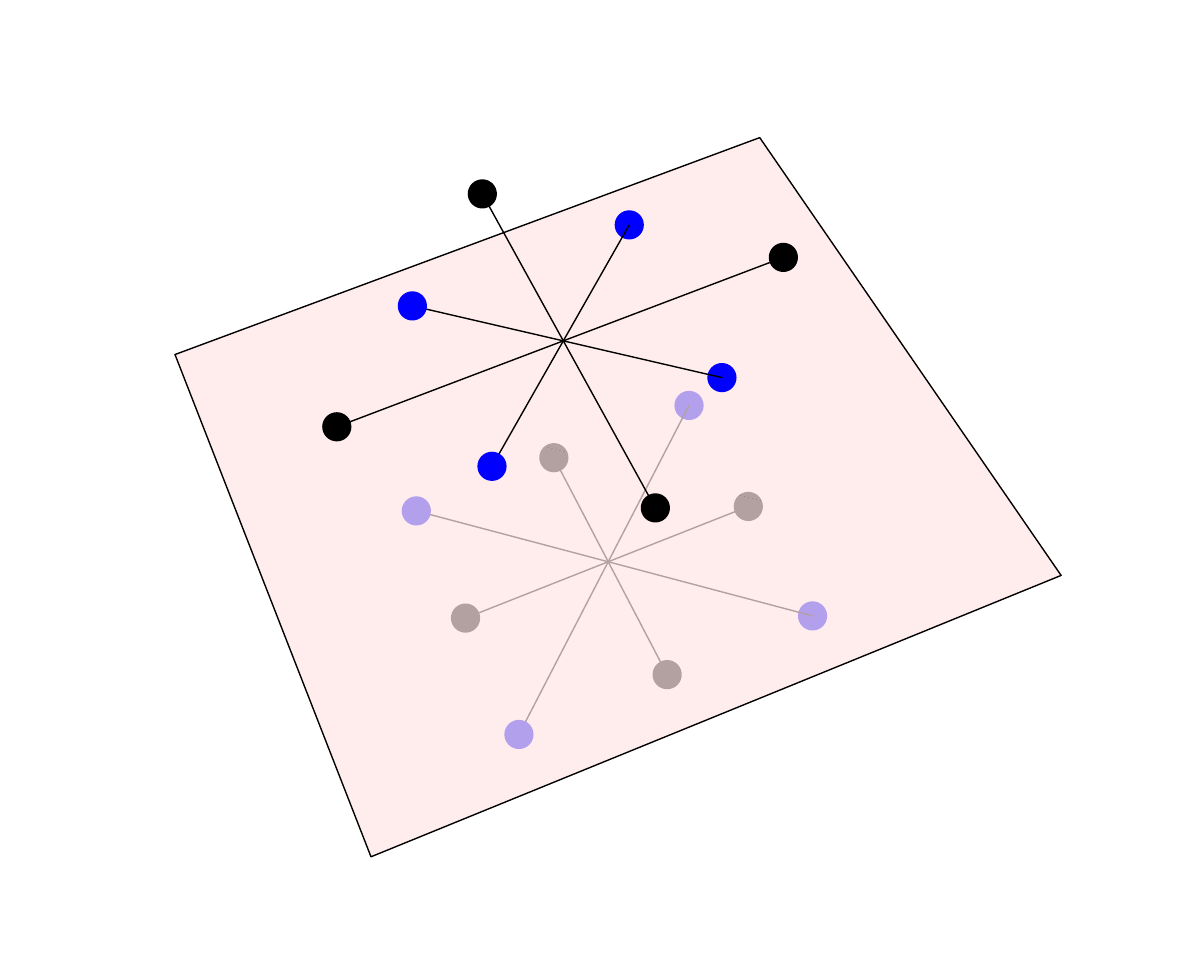}

 \caption{The root system of the group $Sp(4)$ with the roots of the endoscopic group $SO(4)$ highlighted. Reflected in the background is the root system for the Langlands dual group $SO(5)$ with the roots of the subgroup $SO(4)$ highlighted.}
 
\end{figure}

I am grateful to D. Jordan for creating Figures~\ref{orbits} and ~\ref{sp4fig}  which beautifully evoke the mystery of endoscopy.

For their many  insights, I would like to thank S. Gunningham, B. Hannigan-Daley, I. Le, Y. Sakallaridis, T. Schedler, M. Skirvin, X. Zhu, and all the participants in the 2010 Duntroon
Workshop, especially its organizers J. Kamnitzer and C. Mautner.

I would also like to thank B. Conrad, M. Emerton, B. Kra, and M. Strauch for their detailed comments on earlier drafts.

Finally, I gratefully acknowledge the support of NSF grant DMS-0901114 and a Sloan Research Fellowship.


\section{Characters and conjugacy classes}\label{secttrace}

To begin to approach the Fundamental Lemma, let's listen once again
to Langlands~\cite{langlandslettertolang}: 

\begin{quote}
``Nevertheless, it is not the fundamental lemma as such that is critical for the analytic theory of automorphic forms and for the arithmetic of Shimura varieties; it is the stabilized (or stable) trace formula, the reduction of the trace formula itself to the stable trace formula for a group and its endoscopic groups, and the stabilization of the Grothendieck-Lefschetz formula. 
\end{quote}

In this section, we will give a rough impression of the trace formula, and in the next section, 
explain what the term
stable is all about.

%

\subsection{Warmup: finite groups}
To get in the spirit, we begin our discussion with the 
well known character theory of a finite group $G$.
There are many references for the material in this section, for example~\cite{fh}, \cite{serre}.



\begin{defn}
By a {\em representation} of $G$, we will mean a finite-dimensional complex vector space $V$
and a group homomorphism  $\pi:G\to GL(V)$.
\end{defn}
Equivalently, we can form the group algebra $\BC[G] = \{\varphi:G\to \BC\}$ equipped with convolution
$$
\xymatrix{
(\varphi_1 * \varphi_2) (g) = \sum_{g_1 g_2 = g} \varphi_1(g_1) \varphi_2(g_2), & \varphi_1, \varphi_2\in \BC[G],
}$$
and consider finite-dimensional $\BC[G]$-modules.

\begin{example}
(1) Trivial representation: take $V= \BC$ with the trivial action.

(2) Regular representation: take $V=\BC[G]$ with the action of left-translation.

\end{example}

%
%
%

\begin{defn}
The {\em character} of a representation $\pi:G\to GL(V)$ is the function
$$
\xymatrix{
\chi_\pi:G\to \BC 
&
\chi_\pi(g) = \Tr(\pi(g))
}$$
\end{defn}

\begin{defn}
Consider the action of $G$ on itself by conjugation. 

We denote the resulting quotient set by $G\adquot G$ and refer to it as the {\em adjoint quotient}.
Its elements are conjugacy classes $[g] \subset G$.

A {\em class function} on $G$ is a function $f:G\adquot G\to \BC$, or equivalently a conjugation-invariant function
$f:G\to \BC$. 
We denote the ring of class functions by $\BC[G\adquot G]$.

\end{defn}

\begin{lemma}
(1) Each character $\chi_\pi$ is a class function.

(2) Compatibility with direct sums: $\chi_{\pi_1\oplus\pi_2} = \chi_{\pi_1} + \chi_{\pi_2}$.

(3)  Compatibility with tensor products: $\chi_{\pi_1\otimes\pi_2} = \chi_{\pi_1}  \chi_{\pi_2}$.

(4) Trivial representation: $\chi_{triv}(g) = 1$, for all $g$.

(5) Regular representation:
$$
\chi_{reg}(g) = 
\left\{
\begin{array}{cl}
|G|, &   g = e \\
0, & g\not = e
\end{array}\right.
$$
\end{lemma}

Class functions have a natural Hermitian inner product
$$
\xymatrix{
\langle \alpha, \beta\rangle = \frac{1}{|G|}\sum_{g\in G} \ol{\alpha(g)} \beta(g),
& \alpha, \beta\in \BC[G\adquot G].
}$$

\begin{prop}
The characters of irreducible representations of $G$ form an orthonormal basis of the class functions $\BC[G\adquot G]$.
\end{prop}

Thus we have two canonical bases of class functions. On the one hand,
there is the {\em geometric basis} of characteristic functions 
$$
\CO_{[g]}: G\adquot G \longrightarrow \BC
\qquad
\CO_{[g]}(h) = 
\left\{
\begin{array}{cl}
1, &   h \in [g] \\
0, & \text{else}
\end{array}\right.
$$ 
of conjugacy classes $[g]\subset G$. 
These are pairwise orthogonal though not orthonormal since
$\langle \CO_{[g]}, \CO_{[g]} \rangle$ is the volume of the conjugacy class $[g] \subset G$.
On the other hand, there is the {\em spectral basis} of  characters 
$$
\chi_{\iota}: G\adquot G \longrightarrow \BC
\qquad
\chi_{\iota}(h) = \Tr(\pi_\iota(h))
$$ 
of irreducible representations $\pi_\iota$ of $G$.
The geometric basis is something one has on any finite set (though the volumes contain extra information). The spectral basis 
is a reflection of the group structure of the original set $G$.


\begin{remark}
One uses the term {\em spectral} with the following analogy in mind. Given a diagonalizable operator $A$ on a complex vector space $V$, the traditional spectral decomposition
of $V$ into eigenspaces can be interpreted as the decomposition of $V$ into irreducible modules for the algebra $\BC[A]$.
\end{remark}

Given an arbitrary representation $(\pi, V)$, we can expand its character $\chi_\pi$ in the two natural bases to obtain an
identity of class functions. Though $V$ might  be completely mysterious, it nevertheless admits an expansion
into irreducible representations
$$
V \simeq \bigoplus_{\iota\in I} V_\iota^{\oplus m_\iota(\pi)}
$$ 
where $I$ denotes the set of irreducible representations.
Thus we obtain an identity of class functions
$$
\sum_{[g]\in G\adquot G}
 \Tr(\pi(g)) \CO_{[g]} = \sum_{\iota \in I} m_\iota(\pi) \chi_{\iota}
$$

The left hand side is geometric and easy to calculate. The right hand side is spectral both mathematically speaking
and in the sense that like a ghost we know it exists though we may not be able to see it clearly. The formula gives us a starting point to understand the right hand side in terms
of the left hand side.
%


\begin{remark}
It is very useful to view the above character formula as an identity of distributions. Namely, given any test function $\varphi:G\to \BC$,
we can write $\varphi= \sum_{g\in G} \varphi(g) \delta_g$, where $\delta_g:G\to \BC$ is the characteristic function of the group element $g$.
Then since everything in sight is linear, we can evaluate the character formula on $\varphi$. This is very natural from the perspective
of representations as modules over the group algebra $\BC[G]$.
\end{remark}

In general, it is difficult to construct representations. Outside of the trivial and regular representations, the only others
that appear immediately from the group structure of $G$ are induced representations.

\begin{defn}
Fix a subgroup $\Gamma \subset G$.
%


For  a representation $\pi:\Gamma\to GL(W)$,
 the corresponding {\em induced representation}
$\pi_{ind}:G\to GL(V)$ is defined by
$$
\xymatrix{
V = \{f: G \to W | f(\gamma x) = \pi(\gamma) f(x) \}
&
\pi_{ind}(g) f(x) = f(xg)
}
$$
We denote by $\chi_{ind}$ the character of $\pi_{ind}$.
\end{defn}

Calculating the character $\chi_{ind}$ is a particularly simple but salient calculation from Mackey theory. 
Since we have no other starting point, we will focus on the induction of the trivial representation. (Exercise: the induction
of the regular representation of $\Gamma$ is the regular representation of $G$.) When we start with the
trivial representation, the induced
representation 
$$V=\{f:\Gamma\bs G \to \BC\}
$$ is simply the vector space of functions on the finite set $\Gamma\bs G$.
It has a natural basis given by the characteristic functions of the cosets.
Thus every element of $G$, or more generally, the group algebra $\BC[G]$, acts on the vector space $V$ by a matrix
with an entry for each pair of cosets $x, y\in \Gamma\bs G$.

\begin{lemma} 
An element  $\varphi:G\to \BC$ of the group algebra $\BC[G]$  acts on the induced representation $V=\{f:\Gamma\bs G \to \BC\}$
by the matrix
$$
K_\varphi(x, y) = \sum_{\gamma\in \Gamma}  \varphi(x^{-1} \gamma y),
\qquad
x, y\in \Gamma\bs G.
$$

\end{lemma}

Now to calculate the character $\chi_{ind}$, we need only take the traces of the above matrices,
or in other words, the sum of their entries when $x=y \in \Gamma\bs G$.

\begin{corollary}
The character $\chi_{ind}$ is given by the formula
$$
\chi_{ind} (\varphi) = \sum_{\gamma\in \Gamma\adquot\Gamma} a_\gamma  \int_{[\gamma] \subset G} \varphi,
$$
where $a_\gamma$ denotes the volume, or number of elements, of the quotient of centralizers $\Gamma_\gamma\bs G_\gamma$,
and the integral denotes the sum 
$$
\int_{[\gamma] \subset G} \varphi   = \sum_{x\in G_\gamma \bs G} \varphi (x^{-1} \gamma x)
$$
 over the $G$-conjugacy class of $\gamma$.
\end{corollary}

\begin{remark}
Suppose we equip the quotients $G\adquot G$, $\Gamma\adquot \Gamma$ with the natural quotient measures, and
let $p:\Gamma/\Gamma \to G/G$ denote the natural projection. Then the above corollary can be concisely rephrased that  $\chi_{ind}$ is 
the pushforward along $p$ of the quotient measure on $\Gamma\adquot \Gamma$. 
\end{remark}

\begin{defn} For $\gamma\in \Gamma$,
the distribution on $G$ given on a test function $\varphi:G\to \BC$ by the integral
over the conjugacy class
$$
 \CO_\gamma(\varphi) 
= \int_{[\gamma] \subset G} \varphi
  = \sum_{x\in G_\gamma \bs G} \varphi (x^{-1} \gamma x)
 $$
is called an {\em orbital integral}. 
\end{defn}

We have arrived at  the {\em Frobenius character formula} for an induced representation
\begin{equation}\label{finite trace formula}
\sum_{\gamma\in \Gamma\adquot \Gamma} a_\gamma  \CO_\gamma(\varphi) 
= \sum_{\iota \in I} m_\iota(\pi_{ind}) \chi_{\iota}(\varphi)
\end{equation}
This is the most naive form of the Arthur-Selberg trace formula. Observe that the right hand side remains mysterious,
but the left hand side is now a concrete geometric expression involving volumes and orbital integrals.




\subsection{Poisson Summation Formula}

Now let us leave finite groups behind, and consider generalizations of the Frobenius character formula \eqref{finite trace formula}.
We will begin by sacrificing complicated group theory and restrict to the simplest commutative Lie group $G = \BR$.

A deceptive advantage of working with a commutative group is that we can explicitly calculate its spectrum.

\begin{lemma}
The irreducible representations of $\BR$ are the characters
$$
\chi_\lambda:\BR\to \BC^\times
\qquad
\chi_\lambda(x) = \exp(2\pi i \lambda x)
$$
with $\lambda\in \BC$.
In particular, the irreducible unitary representations are the characters $\chi_\lambda$, with $\lambda\in \BR$.
\end{lemma}

Now let us consider the analogue of the Frobenius character formula \ref{finite trace formula} for the group $G= \BR$.
In order for the formula (not to mention our derivation of it) to make immediate sense, we should restrict to a subgroup $\Gamma \subset G$
which is discrete with compact quotient  $\Gamma\bs G$. Thus we are led to the subgroup of integers $\Gamma = \BZ$ with quotient
the circle $S^1 \simeq \BZ\bs \BR$. 

Let us calculate the various terms in the formula \ref{finite trace formula} for the induced
Hilbert representation  $L^2(S^1)$ of square-integrable complex-valued functions. For the geometric side, since $\BR$ is commutative, the conjugacy class of $n\in \BZ$ is simply $n$ itself
with volume $1$.
For the spectral side, Fourier analysis confirms that the representation $L^2(S^1)$ is a Hilbert space direct sum
of the irreducible characters $\chi_\lambda$, with $\lambda\in 2\pi i \BZ$. Furthermore, a compactly supported test function
$\varphi\in C_c^\infty(\BR)$ acts on the summand $\chi_\lambda$ by multiplication  by its Fourier transform
$$
\wh \varphi(\lambda) = \int_{\BR} \varphi(x)\chi_\lambda(x) dx.
$$

\begin{thm}[Poisson Summation Formula]
For a test function $\varphi\in C_c^\infty(\BR)$, one has the equality
$$
\sum_{n\in \BZ} \varphi(n) = 
\sum_{\lambda \in \BZ} 
\wh \varphi(\lambda) 
$$
\end{thm}

With this success in hand, one could seek other commutative groups and attempt a similar analysis.
From the vantage point of number theory, number fields provide a natural source of locally compact
commutative groups.

By definition, a number field $F$ is finite extension of the rational numbers $\BQ$.
There is a deep and
 pervasive analogy between number fields and the function fields $k(X)$ of  algebraic curves $X$. For example, the fundamental example of a number field is $\BQ$, and by definition, all others are finite extensions of it. The fundamental example of an algebraic curve is 
the projective line $\BP^1$, and all other algebraic curves are finite covers of it. The history of the analogy is long with many refinements by celebrated mathematicians (Dedekind, Artin, Artin-Whaples, Weil,...).
 As we will recount below, one of the most intriguing
aspects of the (currently known) proof of the Fundamental Lemma is its essential use of function fields and the extensive analogy between them and number fields.

Throughout what follows, we will need the number field analogue of the  most basic construction of Calculus: the
Taylor series expansion of a function around a point.
Given a curve $X$ and a (rational)  point $x\in X$,  we can choose a local coordinate $t$ with a simple zero at $x$. Then for any non-zero rational function $f\in k(X)$, we have its Laurent series expansion 
$$
\sum_{i = j}^\infty a_i t^i \in k((t)),\text{
with $a_j \not = 0$.}
$$
Since rational functions are locally determined, this provides an embedding of fields $k(X)\subset k((t))$.
The embedding realizes $k((t))$ as the completion of $k(X)$ with respect to the valuation $v_x(f) = k$.

Let us illustrate the form this takes for number fields with the fundamental example of the rational numbers $\BQ$.
The local  expansion of an element of $\BQ$ should take values in a completion of $\BQ$.
Ostrowski's Theorem confirms that the completions are precisely the $p$-adic
numbers $\BQ_p$, for all primes $p$, along with the real numbers $\BR$. 
The real numbers are of course complete
with respect to the usual Euclidean absolute value. 
The $p$-adic numbers are complete with respect to the absolute value $|f|_p = p^{-k}$, where $f= p^{k} a/b$, with $(a, p) = (b, p) =1$.
It satisfies  the non-Archimedean property $|f+g|_p \leq \max\{|f|_p,|g|_p\}$, and so the compact unit ball of $p$-adic integers 
$$
\BZ_p = \{ f \in \BQ_p | |f|_p \leq 1\} \subset \BQ_p
$$ 
is in fact a subring.

It is an elementary but immensely useful idea 
to keep track of all of the local expansions of a rational number
at the same time. Observe that for a rational function on a curve, the points where it has a pole are finite
in number. Similarly, only finitely many primes divide the denominator of a rational number. This leads one to form the ring of adeles
$$
\BA_\BQ = \prod^{rest}_{p \text{ prime}} \BQ_p \times \BR
$$
where the superscript ``rest" denotes that we take the restricted product  where all but finitely many
terms lie in the compact unit ball of $p$-adic integers $\BZ_p$.
The simultaneous local expansion of a rational number provides an embedding of rings
$
\BQ \subset \BA_\BQ
$
with discrete image.
%
%

Let us justify the above somewhat technical construction with a well known result of number theory.
To solve an equation in $\BQ$, it is clearly necessary to provide solutions in 
 $\BQ_p$, for all primes $p$, and also $\BR$.
The Hasse principle asserts that to find solutions in $\BQ$,  one should start with such a solution in the adeles $\BA_\BQ$,
or in other words, a collection of possibly unrelated solutions, and attempt to glue them together. Here is an example
of the success of this approach.

\begin{thm}[Hasse-Minkowski] Given a quadratic form 
$$ 
\xymatrix{
Q(x_1, \ldots, x_n) = \sum_{i \leq j} a_{ij} x_i x_j,
& 
a_{ij} \in \BQ,
}
$$
the equation 
$$
Q(x_1, \ldots, x_n) = 0
$$ 
has a solution in the rational numbers $\BQ$ 
if and only if it has a solution in the adeles $\BA_\BQ$, or equivalently,  solutions in the $p$-adic numbers $\BQ_p$, for all primes $p$, and the real numbers $\BR$.
\end{thm}

A similar constructions of adeles make sense for arbitrary number fields $F$.
The completions of $F$ will be finite extensions  of the completions of $\BQ$, so finite extensions $F_\fp$ of the $p$-adic
numbers $\BQ_p$, along with possibly
 the real numbers $\BR$ or complex numbers $\BC$. 
 The former are non-Archimedean so the compact unit balls of integers $\CO_{\fp} \subset F_\fp$ are in fact subrings.
 One similarly forms the ring of adeles 
 $$
\BA_F = \prod^{rest}_{\fp} F_\fp \times \BR^r \times \BC^c
$$
where $\fp$ runs over all non-Archimedean completions of $F$,
and 
the superscript ``rest" denotes that we take the restricted product  where all but finitely many
terms lie in the compact unit balls $\CO_p$.
The simultaneous local expansion of elements  provides an embedding 
$
F \subset \BA_F
$
with discrete image.

In his celebrated thesis, Tate generalized the Poisson Summation Formula to the pair of the locally compact group $\BA_F$
and its discrete subgroup $F$. 
The resulting formula is an exact
analogue of the classical Poisson Summation Formula
$$
\sum_{x\in F} \varphi(x) = 
\sum_{\lambda \in F} 
\wh \varphi(\lambda) 
$$
This is an essential part of Tate's interpretation of Class Field Theory in terms of harmonic analysis.
One could approach all that follows as an attempt to explore the generalization of Class Field Theory to a noncommutative setting.


\subsection{Arthur-Selberg Trace Formula}

The Arthur-Selberg Trace Formula is a vast generalization of the Frobenius character formula for finite groups
and the Poisson summation formula for number fields. 

The starting point is an algebraic group $G$ defined over a number field $F$.
One can always realize $G$ as a subgroup of $GL(n)$ defined by polynomial equations with coefficients in $F$.
We will be most interested in {\em reductive} $G$, which means that we can realize $G$ as a subgroup of $GL(n)$
 preserved by the transpose of matrices, or equivalently, that the unipotent radical of $G$ is trivial. 
 Without further comment, we will also assume that $G$ is connected in the sense that $G$ is not the union of two proper subvarieties.
 Of course,  $GL(n)$ itself is a fundamental example of a reductive algebraic group. For many important questions, the reader would lose nothing considering only $GL(n)$. But as we shall see, the role of the Fundamental Lemma is to help us compare different groups,
and in particular, reduce questions about complicated groups to simpler ones.

Suppose we are given an algebraic group $G$ defined over a number field $F$. 
Then it makes sense to consider the solutions $G(R)$ to the equations defining $G$ in any commutative ring $R$ containing $F$.
Such solutions are called the {\em $R$-points of $G$} and form a group in the traditional sense of being a set equipped with a group law.

Less naively, although more trivially, for a field $K$ containing $F$, we can also regard the coefficients of the equations defining $G$ as elements of $K$. Hence
we can consider $G$ as an algebraic group defined over $K$. To keep things straight, we will write $G_K$ to denote $G$
thought of as an algebraic group defined over $K$. We will refer to $G_K$ as the {\em base change of $G$}
since all we have done is change the base field.

The only difference between the base change $G_K$ and  the original group $G$ is that
we are only allowed to form the $R$-points $G_K(R)$ of the base change for commutative rings $R$ containing $K$.
Experience tells us that over algebraically closed fields, there is little
difference between equations and their solutions. In practice, this is true for algebraic groups: for an algebraic closure $\ol F$,
one can go back and forth between the $\ol F$-points $G(\ol F)$ and the base change $G_{\ol F}$.

Here is the most important class of reductive groups.

\begin{defn}
(1) A {\em torus} $T$ defined over $F$ is an algebraic group defined over $F$ such that the base change $T_{\ol F}$ is isomorphic to
the product $GL(1)^k$, for some $k$.

(2) A {torus} $T$ defined over $F$ is said to be {\em split} if it is isomorphic to the product $GL(1)^k$, for some $k$,
without any base change necessary.

(3)  A {torus} $T$ defined over $F$ is said to be {\em anisotropic}, or synonymously {\em elliptic}, if all of its characters are trivial
$$
\Hom_F(T, GL(1)) = \langle 0\rangle
$$
where $\Hom_F$ denotes homomorphisms of algebraic groups defined over $F$.
\end{defn}

\begin{example}\label{ex: tori} There are two types of one-dimensional tori defined over the real numbers $\BR$. 
There is the split torus 
$$
GL(1) \simeq \{ab = 1\}
$$
with $\BR$-points $\BR^\times$, and the anisotropic torus 
$$
SO(2) = \{ g\in GL(2) | g^{-1} =  g^\tau \} \simeq \{x^2 + y^2 = 1\}
$$
with $\BR$-points the circle $S^1$. 
Over the complex numbers $\BC$, the equations $ab = 1$ and $x^2 + y^2 =1$ become equivalent via the transformation
$a = x + iy$, $b = x - iy$.
\end{example}

It is often best to think of an algebraic group $G$ defined over $F$ as comprising roughly two pieces of structure:

\begin{enumerate}
\item the base change $G_{\ol F}$ or group of $\ol F$-points $G(\ol F)$ for an algebraic closure $\ol F$, 
\item 
the Galois descent data needed to recover the original equations of $G$ itself. 
\end{enumerate}

To understand the result of the first step, we recall the following definition.

\begin{defn}  \label{def split}

(1) A torus $T\subset G$ is said to be  {\em maximal} if it is not a proper subgroup of another torus in $G$.

(2) A reductive algebraic group $G$ is said to be {\em split} if it contains a maximal torus which is split.
\end{defn}

\begin{prop}
Let $G$ be a reductive algebraic group defined over a number field $F$. Then there is a unique split reductive algebraic group $G^{spl}$
defined over $\BQ$ such that 
$$
G_{\ol F} \simeq G^{spl}_{\ol F},
\quad
\text{ and in particular }
\quad 
G(\ol F) \simeq G^{spl}(\ol F).
$$
In other words, over algebraically closed fields, all reductive groups are split.
\end{prop}

There is a highly developed structure theory of reductive algebraic groups, but the subject is truly example oriented. 
There is the well known Cartan classification of split reductive algebraic groups.

\begin{example}[Split classical groups] 
\label{classical groups}
There are four series of automorphism groups of familiar linear geometry. 
\begin{enumerate}
\item[($A_n$)] The special linear group 
$$
SL({n+1}) = \{A \in GL(n+1) | \det(A) = 1\}.
$$

\item[($B_n$)] The odd special orthogonal group 
$$
SO({2n+1}) = \{A \in GL(2n+1) | A^\tau Q_{2n+1} A  =  Q_{2n+1}, \det(A) = 1\},
$$
$$
Q_{2n+1}
 = \left[
\begin{array}{ccc}
0 & 0 & 1\\
 0 & Q_{2n-1} & 0\\
 1 & 0 & 0
\end{array}
\right].
$$
\item[($C_n$)] The symplectic group 
$$
Sp({2n}) = \{A \in GL(2n) | A^\tau \Omega_{2n} A =  \Omega_{2n}, \det(A) = 1\},
$$
$$
\Omega_{2n}
 = \left[
\begin{array}{ccc}
0 & 0 & -1\\
 0 &  \Omega_{2n-2} & 0\\
 1 &   0& 0
\end{array}
\right].
$$
\item[($D_n$)] The even special orthogonal group 
$$
SO({2n}) = \{A \in GL(2n) | A^\tau Q_{2n}A =  Q_{2n}, \det(A) = 1\},
$$
$$
Q_{2n}
 = \left[
\begin{array}{ccc}
0 & 0 & 1\\
 0 & Q_{2n-2} & 0\\
 1 & 0 & 0
\end{array}
\right].
$$

\end{enumerate}

Each is simple in the sense that (it is not a torus and) 
any normal subgroup is  finite, unlike for example $GL(n)$ which has center $GL(1)$ realized as diagonal invertible matrices.

Each is also {split} with maximal torus its diagonal matrices. 
Outside of finitely many exceptional groups, all simple split reductive algebraic groups are isomorphic to one of the finitely many finite central
extensions or finite quotients of the above classical groups.
\end{example}

To illustrate the Galois descent data involved in recovering a reductive group $G$ defined over $F$ from its base change $G_{\ol F}$,
let us consider the simple example of a torus $T$. The base change $T_{\ol F}$ is split and hence isomorphic to $GL(1)^k$, for some $k$. Its characters form a lattice
$$
X^*(T_{\ol F}) = \Hom(T_{\ol F}, GL(1)) \simeq \BZ^k
$$
from which we can recover $T_{\ol F}$ as the spectrum 
$$
T_{\ol F} = \Spec \ol F[X^*(T_{\ol F})].
$$
The Galois group $\Gamma = Gal(\ol F/F)$ naturally acts on 
the character lattice $
X^*(T_{\ol F}) 
$
by a finite group of automorphisms. This induces an action on the ring $ \ol F[X^*(T_{\ol F})]$, and 
we can recover $T$ as the spectrum of the ring of invariants
$$
T = \Spec \ol F[X^*(T_{\ol F})]^\Gamma.
$$
In the extreme cases,  $T$ is split if and only if the Galois action on $X^*(T_{\ol F})$ is trivial, and by definition,
$T$ is anisotropic if and only if
the invariant characters $X^*(T_{\ol F})^\Gamma = \Hom_F(T, GL(1))$ are trivial.

%
%
%
%

There is a large class of non-split groups which contains all tori and is particularly easy to describe by Galois descent.
All groups directly relevant to the Fundamental Lemma will come from this class. We give the definition and a
 favorite example here, but
defer discussion of the Galois descent until Section~\ref{sect: endo}

\begin{defn} \label{def quasisplit}

(1) A {\em Borel subgroup} $B\subset G$ is an algebraic subgroup such that the base change $B_{\ol F}\subset G_{\ol F}$
is a maximal solvable algebraic subgroup. 
%
%
 
(2) A reductive algebraic group $G$ is said to be {\em quasi-split} if it contains a Borel subgroup $B\subset G$.
\end{defn}

\begin{example}
[Unitary groups] Suppose $E/F$ is a separable degree $2$ extension of fields. Then there is a unique nontrivial
involution of $E$ fixing $F$ which one calls conjugation and denotes by $a\mapsto \ol a$, for $a\in E$. The unitary group 
is the matrix group
$$
U({n, J, E/F}) = \{A \in GL(n, E) | \ol A^\tau J A =  J, \det(A) = 1\}
$$
where $\ol J^\tau = J$ is nondegenerate.
We can view $U({n, J,  E/F})$ as a subgroup of $GL(2n, F)$ cut out by  equations defined over $F$,
and hence $U({n, J, E/F})$  is a reductive algebraic group defined over $F$.

If we take $J$ to be antidiagonal, then $U({n, J,  E/F})$ is quasi-split with Borel subgroup its upper-triangular matrices.
But if we take $J$ to be the identity matrix, then 
for example $U({n, J, \BC/\BR})$ is the familiar compact unitary group which is as far from quasi-split as possible.
(All of its connected algebraic subgroups are reductive, and so have trivial unipotent radical.)
\end{example}

%
%
%
%

Suppose we are given an algebraic group $G$ defined over a number field $F$,
and also a commutative ring $R$ with a locally compact topology. Then the group of $R$-points $G(R)$ is a locally compact topological group, and so amenable
to the techniques of harmonic analysis. Most prominently, $G(R)$ admits a bi-invariant Haar measure, and we can study
representations such as 
$L^2(\Gamma\bs G(R))$, for discrete subgroups $\Gamma\subset G(R)$.
As harmonic analysts, our hope is to classify the irreducible representations of $G(R)$, and decompose
complicated representations such as $L^2(\Gamma\bs G(R))$ into irreducibles. Our wildest dreams are
to find precise analogues of the successes  of Fourier analysis where the initial group $G$ is commutative. 

For example, suppose we are given a reductive algebraic group $G$ defined over the integers $\BZ$,
so in particular, the rational numbers $\BQ$. 
Then we can take $K$ to be either the local field of $p$-adic
numbers $\BQ_p$ or real numbers $\BR$. We obtain the $p$-adic groups $G(\BQ_p)$ and the real Lie
group $G(\BR)$. They are locally compact with respective maximal compact subgroups $G(\BZ_p)$
where $\BZ_p\subset \BQ_p$ is the compact
unit ball of $p$-adic integers, and $K(\BR)$, where $K\subset G$ is the fixed points of the involution which takes a matrix
to its inverse transpose.
Finally, we can consider them all simultaneously by forming the locally compact ad\`elic group 
$$
G(\BA_\BQ)
\simeq \prod^{rest}_{p \text{ prime}} G(\BQ_p) \times G(\BR).
$$
It is a fundamental observation that the inclusion $\BQ \subset \BA_\BQ$ induces an inclusion $G(\BQ) \subset G(\BA_\BQ)$ with discrete image, and so the space of automorphic functions 
$
L^2(G(\BQ)\bs G(\BA_\BQ))
$
presents a natural representation of $G(\BA_\BQ)$, and in particular of the $p$-adic groups $G(\BQ_p)$ and 
Lie group $G(\BR)$, 
to approach via harmonic analysis.

In general, for a reductive algebraic group $G$ defined over an arbitrary number field $F$, 
by passing to all of the completions of $F$, we obtain the locally compact $p$-adic groups $G(F_\fp)$ and possibly the Lie
groups $G(\BR)$ and $G(\BC)$, depending on whether $\BR$ and $\BC$ occur as completions.
They are locally compact with respective maximal compact subgroups $G(\CO_\fp) \subset G(F_\fp)$
where $\CO_\fp\subset F_\fp$ is the ring
of integers, $K(\BR)\subset G(\BR)$, where $K\subset G$ is the fixed points of the involution which takes a matrix
to its inverse transpose, and $U(\BC)\subset G(\BC)$, where $U\subset G$ is the fixed points of the involution which takes a matrix to its conjugate inverse transpose. 
 We can form the locally compact ad\`elic group 
 $$
G(\BA_F) \simeq \prod^{rest}_{\fp} G(F_\fp) \times G(\BR)^r \times G(\BC)^c
$$
where $\fp$ runs over all non-Archimedean completions of $F$.
As with the rational numbers, the inclusion $F \subset \BA_F$ induces an inclusion $G(F) \subset G(\BA_F)$ with discrete image,
and so the space of automorphic functions 
$
L^2(G(F)\bs G(\BA_F))
$
presents a natural representation of $G(\BA_F)$
to approach via harmonic analysis.

\begin{example}
The automorphic 
representation $L^2(G(F)\bs G(\BA_F))$ is far less abstract than might initially appear. Rather than recalling the general statement of {strong approximation}, we will focus on the classical case when our number field is the rational numbers $F=\BQ$ and our group is $G=SL(2)$. Inside of the ad\`elic group $G(\BA_\BQ)$, consider the product of maximal compact subgroups
$$
K = 
 \prod_{p \text{ prime}} SL(2, \BZ_p) \times SO(2, \BR).
$$
Then with $\BH \subset \BC$ denoting the open upper halfplane, there is a canonical identification
\begin{eqnarray*}
SL(2, \BQ)\bs SL(2, \BA_\BQ) / K& \simeq & SL(2, \BZ) \bs SL(2, \BR)/ SO(2, \BR) \\
& \simeq  & SL(2, \BZ) \bs \BH,
\end{eqnarray*}
and the latter is the moduli of elliptic curves. By passing to smaller and smaller subgroups of $K$,
we obtain the moduli
of elliptic curves with level structure.
This classical realization opens up the study of the original automorphic representation to the more 
familiar techniques (Laplace-Beltrami operators, Hecke integral operators) of harmonic analysis.
\end{example}

\begin{remark}
It is beyond the scope of this article to explain, but suffice to say, the deepest secrets of the universe are contained in the spectrum
of the automorphic representation $L^2(G(F)\bs G(\BA_F))$. The {\em Langlands correspondence} is a conjectural description of the spectrum in terms of 
representations of the Galois group $Gal(\ol F/F)$. Thanks to the symmetry of the situation, one can turn things around and attempt to understand 
$Gal(\ol F/F)$ in terms of  $L^2(G(F)\bs G(\BA_F))$. When one shows that a Galois representation is automorphic, or in other words, occurs in the 
 spectrum, this leads to many deep structural implications. 
 
Not only in general, but even in specific cases, it is extremely difficult to confirm
that a given Galois representation  is automorphic. Often the only hope is to bootstrap off of the precious few historical successes
by concrete techniques such as  induction and less obviously justified approaches such as prayer.
But the prospect of success is at least supported by {\em Langlands's   functoriality} which conjectures that
whenever there is an obvious relation between Galois representations, there should be a parallel relation between automorphic representations. In particular, there are  often highly surprising relations between automorphic representations for different groups
corresponding to much more prosaic relations of Galois representations. It is in this context that the Fundamental Lemma plays an essential role.
\end{remark}

Now we arrive at the Arthur-Selberg Trace Formula which is the primary tool in the study of automorphic representations. 
For simplicity, let us restrict for the moment to the far more elementary setting where the quotient $G(F)\bs G(\BA_F)$ is compact.
The group algebra $C_c^\infty(G(\BA_F))$ of smooth, compactly supported functions on the ad\`elic group acts on the automorphic representation $L^2(G(F)\bs G(\BA_F))$ by compact operators
$$
R(\varphi)f(g) = \int_{G(\BA_F)}	\varphi (gh)f(h)dh.
$$
It follows that the representation decomposes as a Hilbert space direct sum of irreducible unitary representations
$$
L^2(G(F)\bs G(\BA_F)) \simeq \bigoplus_\iota m_\iota \pi_\iota
$$
We can form the character $
\Tr R
$
as a distribution on $G(\BA_F)$.
The formal analogue of
the Frobenius character formula~\ref{finite trace formula} is an instance of the Selberg Trace Formula.

\begin{thm}[Selberg Trace Formula for compact quotient] Suppose $G(F)\bs G(\BA_F)$ is compact. 
Then for any test function $\varphi \in C_c^\infty(G(\BA_F))$,  we have an identity
$$
\sum_{\gamma\in G(F)\adquot G(F)} a_\gamma \CO_\gamma(\varphi)
=
\sum_{\iota } m_\iota \chi_\iota(\varphi)
$$
where $a_\gamma$ is the volume of the quotient $G_\gamma(F)\bs G_\gamma(\BA_F)$, and 
 the distribution $\CO_\gamma$ is the orbital integral 
$$
 \CO_\gamma(\varphi) 
= \int_{[\gamma] \subset G} \varphi
 $$
over the $G(\BA_F)$-conjugacy class $[\gamma]\subset G(\BA_F)$.
\end{thm}

For modern theory and applications, one needs Arthur's generalizations of the Selberg Trace Formula for very general quotients. The technical details are formidable and Arthur's expositions can not be improved upon. But its formal structure 
and application is the same.
On the geometric side, we have a formal sum involving volumes and explicit orbital integrals in the ad\`elic group.
On the spectral side, we have the character of the automorphic representation expressed as a formal integral over the characters  of irreducible representations. 
The identification
of the two sides gives us a starting point to attempt to understand the spectrum in terms of geometry.

Although there are important and difficult issues in making this formal picture rigorous,
there is an immediately accessible piece of it which can be isolated. On the spectral side, there is the discrete part of the automorphic
spectrum consisting of irreducibles which occur on their own with positive measure. On the geometric side, there are orbital integrals for elements $\gamma \in G(F)$ whose centralizers $G_\gamma$ are   {anisotropic}
tori.
 
The Fundamental Lemma is needed for the comparison of the anisotropic terms of the geometric side of the   trace formula for different groups.
We can leave for another time the thorny complications of other aspects of the trace formula. From hereon,
we can focus on orbital integrals over anisotropic conjugacy classes. Moreover, we can expand
each anisotropic orbital integral around each ad\`elic place to obtain 
\begin{eqnarray}\label{adelic factorization}
\CO_\gamma(\varphi) = \prod_{\fp} \CO_{\gamma}(\varphi_\fp)
\end{eqnarray}
where $\fp$ runs over all completions of $F$, we expand $\gamma$ at each place, $\CO_\gamma(\varphi_\fp)$ denotes the orbital integral along the conjugacy class $[\gamma] \subset G(F_\fp)$, and without sacrificing too much, we work with a product test function $\varphi = (\varphi_{\fp})$.
Thus from hereon, leaving global motivations behind, we can focus on orbital integrals over conjugacy classes in local groups.

\section{Eigenvalues versus characteristic polynomials}\label{sectstable}

Our discussion of the previous section is a success if the reader comes away with the impression that outside
of the  formidable technical issues in play, the basic idea of the trace formula is a kind of formal tautology.
The great importance and magical applications of Arthur's generalizations to arbitrary ad\`elic groups are found in comparing trace formulas for
different groups. This is the primary approach to realizing instances of Langlands's functoriality conjectures
on the relation of automorphic forms on different groups. The general strategy is to compare the geometric sides
where traces are expressed in concrete terms, and thus arrive at conclusions about the mysterious spectral sides.
By  instances of Langlands's reciprocity conjectures, the spectral side involves Galois theory, and 
eventually leads to deep implications in number theory.

Now an immediate obstruction arises when one attempts to compare the geometric sides of the trace formulas for different groups. Orbital integrals over conjugacy classes in different groups  have no evident relation with each other. 
Why should we expect conjugacy classes of say symplectic matrices and  orthogonal matrices to have anything to talk about? If we diagonalize
them, their eigenvalues live in completely different places. But here is the key observation that gives one hope: 
{\em the equations describing their eigenvalues
are in fact intimately related}. 
In other words, if we pass to an algebraic closure, where equations and their solutions are more closely tied, then we 
find a systematic relation between conjugacy classes. To explain this further, we will start with some elementary linear algebra, then
build to Langlands's theory of endoscopy, and in the end,  arrive at the Fundamental Lemma.

\subsection{The problem of Jordan canonical form}

Suppose
we consider a field $k$, and a finite-dimensional $k$-vector space $V\simeq k^n$. Given an endomorphism $A\in \End_k(V)\simeq M_{n\times n}(k)$, form the characteristic polynomial 
$$
p_A(t) = \det(t\Id_V - A)  = a_0 + a_1 t+ \cdots + a_{n-1} t^{n-1} + t^n \in k[t].
$$ For simplicity, we will assume that the roots
of $p_A(t)$, or equivalently, the eigenvalues $\lambda_1, \ldots, \lambda_n$ of $A$, are all distinct. 
Of course, if $k$ is not algebraically closed, or more generally, does not contain the roots of $p_A(t)$, we will need to pass to an extension of $k$ to speak concretely of the roots.

Let's review the two ``canonical" ways to view
the endomorphism $A$. On the one hand, we can take the coefficients of $p_A(t)$ and form the companion matrix
$$
C_A = \left[
\begin{array}{ccccc}
 0 & 0 & \cdots & 0 & -a_0\\
  1 & 0 & \cdots & 0 & -a_1\\
   0 & 1 & \cdots & 0 & -a_2\\
   \vdots   &\vdots & \ddots & &\vdots\\
    0 & 0 & \cdots & 1 & -a_{n-1}
\end{array}
\right]
$$
Since we assume that $p_A(t)$ has distinct roots, and hence is equal to its minimal polynomial, $C$ is
the rational normal form of $A$, and hence
$A$ and $C_A$ will be conjugate. We think of this as the naive {\em geometric} form of $A$.

On the other hand,
we can try to find a basis of $V$ in which $A$ is as close to diagonal as possible. If $k$ is algebraically closed, or more generally, contains the eigenvalues of $A$, then we will be able to conjugate $A$ into Jordan canonical form. 
In particular, since we assume that $A$ has distinct eigenvalues, $A$
 will be conjugate to the diagonal matrix
$$
D_A = \left[
\begin{array}{cccc}
 \lambda_1 & 0 & \cdots & 0\\
 0 & \lambda_2  & \cdots & 0\\
 \vdots & \vdots & \ddots & \vdots\\
 0 & 0 & \cdots & \lambda_n\\
\end{array}
\right]
$$ 
We think of this as the sophisticated {\em spectral }  form of $A$. It is worth noting that the most naive ``trace formula" is found in the identity
$$
\Tr(A) = -a_{n-1} = \lambda_1 + \cdots + \lambda_n
$$
which expresses the spectral eigenvalues of $A$  in terms of the geometric sum along the diagonal of $A$.

When $k$ is not algebraically closed, or more specifically, does not contain the eigenvalues of $A$, understanding the structure of $A$ is more difficult. It is always possible to conjugate $A$
into rational normal form, but not necessarily Jordan canonical form.
One natural solution is to fix an algebraic closure $\bar{k}$, and regard $A$ as an endomorphism of 
the extended vector space $\ol{V} = V \otimes_k \bar{k} \simeq \bar k^n$. Then we can find a basis of $\ol V$ for which $A$ is in Jordan canonical form. Equivalently, we can conjugate $A$ into Jordan canonical form by an element of the automorphism group $\Aut_{\bar k}(\ol V) \simeq GL(n, \bar k)$. 
This is particularly satisfying since {Jordan canonical forms of matrices completely characterize their structure}.

\begin{lemma}
If two matrices $A, A'\in M_{n\times n}(k)$ are conjugate by an element of $GL(n,\bar k)$, they  are in fact conjugate by an element of $GL(n, k)$.
\end{lemma}

All of the subtlety of the Fundamental Lemma emanates from the difficulty that when we consider subgroups of $GL(n)$, the above lemma consistently fails. For example, suppose we restrict the automorphism group of our vector space $V$ to be the special linear group $SL(n)$. In other words, we impose that the symmetries of $V$ be not all invertible linear maps, 
but only those preserving volume. Then Jordan canonical form is no longer a complete invariant for the equivalence
classes of matrices.

\begin{example}
Take $k=\BR$. Consider the rotations of the real plane $V = \BR^2$ given by the matrices
$$
A(\theta) =  \left[
\begin{array}{cc}
\cos(\theta) & -\sin(\theta) \\
\sin(\theta)  & \cos(\theta)
\end{array}
\right]
\qquad
A'(\theta) =  \left[
\begin{array}{cc}
\cos(\theta) & \sin(\theta) \\
-\sin(\theta)  & \cos(\theta)
\end{array}
\right].
$$
Observe that both $A(\theta), A'(\theta)$ lie in $SL(2, \BR)$, and they are conjugate by the matrix 
$$
M =  \left[
\begin{array}{cc}
i & 0 \\
0  & -i
\end{array}
\right] \in SL(2, \BC).
$$
Furthermore, when $\theta \not\in \pi\BZ$, there is no element in $SL(2, \BR)$ which conjugates one into the other.
When we view $A(\theta), A'(\theta)$ as endomorphisms of the complex plane $\ol V = \BC^2$, 
they both are conjugate to the diagonal matrix
$$
D(\theta)=  \left[
\begin{array}{cc}
\cos(\theta) + i\sin(\theta) & 0 \\
0  & \cos(\theta) - i\sin(\theta)
\end{array}
\right].
$$
\end{example}

Let us introduce some Lie theory to help us think about the preceding phenomenon.
For simplicity, we will work with a split reductive group $G$ whose derived group 
$G_{der} = [G,G]$ is simply connected. 
For example, the split classical groups of Example~\ref{classical groups} are all simple, hence equal to their derived groups.
The special linear and symplectic groups
are simply-connected, but for the special orthogonal group, one needs to pass to the spin two-fold cover.

Fix a  split maximal torus $T\subset G$, and recall that
the Weyl group of $G$ is the finite group $W=N_T/T$, where $N_T\subset G$ denotes the normalizer
of $T$. All split tori are conjugate by $G(k)$ and the choice of $T\subset G$ is primarily for convenience.
%

To begin, let us recall the generalization of Jordan canonical form. Recall  that to diagonalize matrices with distinct eigenvalues, 
in general, we have to pass to an algebraically closed field $\ol k$.

\begin{defn} For an element $\gamma\in G(k)$, let $G_\gamma\subset G$ denote its centralizer.

(1) The element $\gamma $ is said to be  {\em regular} if  $G_\gamma$ is commutative.

(2) The element $\gamma $ is said to be  {\em semisimple} if  $G_\gamma$ is connected and reductive.

(3) The element $\gamma$ is said to be {\em regular semisimple} if it is regular and semisimple,
or equivalently  $G_\gamma$ is a torus.

(3) The element $\gamma$ is said to be {\em anisotropic} if  $G_\gamma$ is an anisotropic torus.

\end{defn}

\begin{example} Take $k=\BR$ and $G= SL(2)$. Consider the elements
$$
r=  \left[
\begin{array}{cc}
1 & 1 \\
0  & 1
\end{array}
\right]
\qquad
s=  \left[
\begin{array}{cc}
1 & 0 \\
0  & 1
\end{array}
\right]
\qquad
h=  \left[
\begin{array}{cc}
2 & 0 \\
0  & 1/2
\end{array}
\right]
\qquad
a=  \left[
\begin{array}{cc}
0 & -1 \\
1  & 0
\end{array}
\right]
$$
with respective centralizers $T_r= \BZ/2\BZ \times \BA^1$, $T_s= SL(2)$, $T_h= GL(1)$, $T_a= SO(2)$. Thus $r$ is regular, $s$ is semisimple, $h$ and $a$ are regular semisimple, and 
 $a$ is anisotropic.
To see the latter fact, observe that there are no nontrivial homomorphisms  $S^1\to \BR^\times$ of their groups of $\BR$-points.
\end{example}

\begin{remark} Some prefer the phrase {\em strongly regular semisimple}  for an element $\gamma\in G(k)$ whose centralizer $G_\gamma$ is a torus
and not a possibly disconnected commutative reductive group.
When $G_{der}$ is simply-connected, if the centralizer $G_\gamma$ is a  reductive group then it will be connected.
\end{remark}

\begin{remark}
Some might prefer to define anisotropic to be slightly more general. Let us for the moment call  
 a regular semisimple element  $\gamma \in G(F)$ {\em anistropic modulo center} if the quotient $T/Z(G)$ of the centralizer 
 $T = {G}_\gamma$ by the center $Z(G)$ is anisotropic. 
 
 A group with split center such as $GL(n)$ will not have anisotropic elements, but will have elements anisotropic modulo center. A regular semisimple element $\gamma\in GL(n)$ will be anisotropic modulo center if and only if its characteristic polynomial is 
 irreducible (and separable).
 
For the Fundamental Lemma, we will be able to focus on anisotropic elements. Somewhat surprisingly, it is not needed for $GL(n)$ where there is no elliptic endoscopy.
\end{remark}


The following justifies the idea that semisimple elements are ``diagonalizable" and regular semisimple elements
are ``diagonalizable with distinct eigenvalues".

\begin{prop}\label{prop diag}
(1) Every semisimple element of $G(\ol k)$ can be conjugated into $T(\ol k)$.

(2) Two semisimple elements of $T(\ol k)$ are conjugate in $G(\ol k)$ if and only if they are conjugate by the Weyl group $W_G$.

(3) An element of $T(\ol k)$ is  regular semisimple if and only the Weyl group $W_G$ acts on it with trivial stabilizer.
\end{prop}

Second, let us generalize the notion of characteristic polynomial. Recall that the coefficients of the characteristic polynomial
are precisely the conjugation invariant polynomial functions on matrices.

\begin{thm}[Chevalley Restriction Theorem]
The $G$-conjugation invariant polynomial functions on $G$ are isomorphic to the $W$-invariant polynomial functions on $T$.
More precisely, restriction along the inclusion $T\subset G$ induces an isomorphism
$$
\xymatrix{
k[G]^G \ar[r]^-\sim & k[T]^W. 
}$$
\end{thm}

Passing from polynomial functions to algebraic varieties, we obtain the $G$-invariant Chevalley morphism
$$
\xymatrix{
\chi:G \ar[r] & T/W   = \Spec  k[T]^W. 
}
$$
It assigns to a group element its ``unordered set of eigenvalues", or in other words its characteristic polynomial.

Finally, let us mention the generalization of rational canonical form for split reductive groups.
Recall that a {\em pinning} of a split reductive group $G$ consists of a Borel subgroup $B\subset G$, split maximal torus $T\subset B$, and basis vectors for the resulting simple positive root spaces. The main consequence of a pinning is that only central  conjugations preserve it, and so it does away with ambiguities coming from inner automorphisms. (For slightly more discussion,
including the example of $G= SL(n)$, see Section~\ref{sect: endo} below where we discuss root data.)

\begin{thm}[Steinberg section] Given a pinning of the split reductive group $G$ with split maximal torus $T\subset G$,
there is a canonical section
$$
\xymatrix{
\sigma: T/W  \ar[r] & G
}
$$
to the Chevalley morphism $\chi$. In other words, $\chi\circ\sigma$ is the identity.
\end{thm}

Thus to each ``unordered set of eigenvalues", we can assign a group element with those eigenvalues. 

With the above results in hand, we can now introduce the notion of stable conjugacy. Recall that
given  $\gamma\in G(k)$, we denote by $[\gamma] \subset G(k)$ the conjugacy class through $\gamma$.

\begin{defn} \label{def stable conj}
Let $G$ be a simply-connected reductive algebraic group defined over a field $k$. 

We say two regular semisimple elements $\gamma, \gamma'\in G(k)$ are {\em stable conjugate} and write $\gamma {\sim_{st}}\gamma'$ if they satisfy one of 
the following equivalent conditions:

(1) $\gamma$ and $\gamma'$ are conjugate by an element of $G(\ol k)$,

(2) $\gamma$ and $\gamma'$ share the same characteristic polynomial
$
\chi(\gamma) = \chi(\gamma').
$

Given $\gamma\in G(k)$, the {\em stable conjugacy class} $[\gamma]_{st} \subset G(k)$ through $\gamma$ consists of  $\gamma'\in G(k)$
stably conjugate to $\gamma$.

\end{defn}

\begin{remark}
Experts in algebraic groups over fields of finite characteristic  will note the useful fact that in Proposition~\ref{prop diag} and Definition~\ref{def stable conj}, one need only go to the separable closure $k_s$.
\end{remark}

Recall that  the geometric side of the trace formula leads to orbital integrals over conjugacy classes of regular semisimple elements in groups over local fields.
The theory of canonical forms for elements
 is intricate, and conjugacy classes are not characterized by their Jordan canonical forms. The complication sketched above
 for  $SL(2, \BR)$   is quite ubiquitous, and one will also encounter it for the classical groups of Example~\ref{classical groups}. Any hope to understand
conjugacy classes in concise terms must involve passage to the stable conjugacy classes found over the algebraic closure. In simpler terms, we must convert constructions  depending on eigenvalues, such as orbital integrals, into constructions depending on
characteristic polynomials.


\subsection{Fourier theory on conjugacy classes}

Imbued with a proper fear of the intricacy of conjugacy classes over non-algebraically closed fields, one dreams that the geometric side of the trace formula
could be rewritten in terms of stable conjugacy classes which are independent of the field.

One might not expect that the set of conjugacy classes in a given stable
conjugacy class would be highly structured.
But it turns out there is extra symmetry governing the situation.

%
%
%

To simplify the discussion, it will be useful to make the standing assumption that the reductive group $G$
is simply connected, or more generally, its derived group $G_{der}=[G,G]$ is simply-connected.

\begin{prop}Let $G$  be a  reductive algebraic group defined over a local field $F$.

Let $\gamma_0\in G(F)$ be a regular semisimple element with centralizer the torus $T = G_{\gamma_0}$. 

Then the set  of conjugacy classes in the stable
conjugacy class $[\gamma_0]_{st}$ is naturally a finite abelian group given by the kernel
of the Galois cohomology map
$$
\xymatrix{
\CA_{\gamma_0}\simeq \ker\{H^1(F, T) \ar[r] & H^1(F, G)\}.
}$$
In particular, it only depends on the centralizing torus $T$ as a subgroup of $G$.

\end{prop}

\begin{remark} (1) When $G$ is simply connected, $H^1(F, G)$ is trivial.

(2)  One can view the Galois cohomology $H^1(F, T)$ as parameterizing principal $T$-bundles over $\Spec F$.
Since $T$ is abelian, this is naturally a group. Under the isomorphism of the proposition, the trivial bundle corresponds to $\gamma_0$.

(3) One can view the quotient $[\gamma_0]_{st}/G(F)$
of the stable conjugacy class by conjugation as a discrete collection
of classifiying spaces for stabilizers. Each of the classifiying spaces is noncanonically isomorphic to the classifying space of $T = G_{\gamma_0}$. The possible isomorphisms form a principal $T$-bundle giving the corresponding class in $H^1(F, T)$.
\end{remark}

Now suppose we have a $G(F)$-invariant distribution on the stable conjugacy class containing $\gamma_0$.
In other words, we have a distribution of the form
$$
\delta(\varphi) = \sum_{\gamma \in \CA_{\gamma_0}} c_\gamma \CO_\gamma(\varphi)
$$
where as usual $\CO_\gamma(\varphi)$ denotes the orbital integral over the $G(F)$-conjugacy class through $\gamma$.
In what sense could we demand the distribution $\delta$ be invariant along the entire stable conjugacy class?
Requiring the coefficients $c_\gamma$ are all equal is a lot to ask for, but there is a  reasonable 
generalization presented by Fourier theory.

Consider the Pontryagin dual group of characters 
$$
\CA^D_{\gamma_0} = \Hom(\CA_{\gamma_0},\BC^\times).
$$

\begin{defn}
Let $G$ be a reductive algebraic group defined over a local field $F$. Let $\gamma_0\in G(F)$
be an anisotropic element.

Given $\kappa\in \CA^D_{\gamma_0}$, the {\em $\kappa$-orbital integral} through $\gamma_0$ is the distribution
$$
\CO^\kappa_{\gamma_0}(\varphi)
=
\sum_{\gamma \in \CA_{\gamma_0}} 
\kappa(\gamma) \CO_\gamma(\varphi)
$$
In particular, when $\kappa=e \in \CA^D_{\gamma_0}$ is the trivial character,  the stable orbital integral is the distribution
$$
\CS\CO_{\gamma_0}(\varphi) =\CO^e_{\gamma_0}(\varphi).
$$
\end{defn}

\begin{remark}
Observe that the stable orbital integral $\CS\CO_{\gamma_0}(\varphi)$ is independent of the choice of base point $\gamma_0$ in the stable conjugacy class.
Thus it is truly associated to the characteristic polynomoial $\chi(\gamma_0)$.

On the other hand, the dependence of the $\kappa$-orbital integral $\CO^\kappa_{\gamma_0}(\varphi)$ on the base point $\gamma_0$ is modest but nontrivial. 
If one chooses some other $\gamma'_0 \in \CA_{\gamma_0}$, the resulting expression will scale by $\kappa(\gamma'_0)$.
For groups with simply-connected derived groups, there is the base point, which is canonical up to a choice of pinning, given by the image of the Steinberg section
$\sigma(\chi(\gamma_0))$.
\end{remark}

Now by Fourier theory, we can write our original distribution $\delta$ as a finite sum
$$
\delta(\varphi) = \sum_{\kappa \in \CA^D_{\gamma_0}} c_\kappa \CO^\kappa_{\gamma_0}(\varphi).
$$
Hence while $\delta$ might not have been stable, it can always be written as a linear combination of  distributions which vary along the stable conjugacy class by a character.

%
%
%
%
%
%
 
 Now to proceed any further, we must understand the character group 
 $$
 \CA^D_{\gamma_0} = \Hom(\CA_{\gamma_0},\BC^\times).
 $$
 A closer examination
 of the possible characters  will reveal the possibility of a deep reinterpretation of the $\kappa$-orbital integrals.
 
Suppose the local field
 $F$ is non-Archimedean, and fix a torus $T$ defined over $F$. Recall that we can think of $T$ as the information of a split torus $T_{\ol F} \simeq GL(1)^k$ over
 the algebraic closure $\ol F$,
 together with the Galois descent data needed to recover the original equations cutting out $T$. The descent
 is captured by the finite action of the Galois group $\Gamma = Gal(\ol F/F)$ on the cocharacter lattice 
 $$
 X_*(T_{\ol F}) = \Hom(GL(1), T_{\ol F}) \simeq \BZ^k.
 $$

 Consider the dual complex torus 
 $$
 T^\vee = \Spec \BC[X_*(T_{\ol F})]
 \simeq GL(1)^k
 $$
 whose monomial functions are the cocharacter lattice. The  $\Gamma$-action on $X_*(T_{\ol F})$ induces a  corresponding $\Gamma$-action on $T^\vee$.
 
 \begin{prop}[Local Tate-Nakayama duality] Assume $F$ is a non-Archimedean local field.
 There is a canonical identification of abelian groups
 $$
H^1(F, T)^D \simeq \pi_0((T^\vee)^\Gamma)
$$
between the Pontryagin dual of  the Galois cohomology of $T$,
and the component group of the $\Gamma$-invariants in the dual torus $T^\vee$.
 \end{prop}

\begin{remark}  When $G$ is simply connected, $H^1(F, G)$ is trivial, and so we have calculated $\CA_{\gamma_0}^D$.

When $G$ is not simply-connected,  elements of $\pi_0((T^\vee)^\Gamma)$ nonetheless restrict to characters of 
 $\CA_{\gamma_0}^D$. It is an exercise to relate the kernel of this restriction to $\pi_1(G)$.
\end{remark}

Thus a regular semisimple element $\gamma_0\in G(F)$ provides a centralizing torus $T= G_{\gamma_0}$ which in turn determines a Galois action on the dual torus $T^\vee$.  To each element $\kappa \in (T^\vee)^\Gamma$ in the Galois-fixed locus, we can associate the $\kappa$-orbital integral $\CO^\kappa_{\gamma_0}(\varphi)$ defined by the image of $\kappa$ in the component group $\pi_0((T^\vee)^\Gamma)$.

%
%
%
%


\subsection{Endoscopic groups and the Fundamental Lemma} \label{sect: endo}

We have reached a pivotal point in our discussion. Let's step back for a moment and take measure of its successes and shortcomings. 

Given a number field $F$,
and a reductive algebraic group $G$ defined over $F$, we aim to 
understand the automorphic representation $L^2(G(F)\bs G(\BA_F))$.
Our main tool is the Arthur-Selberg Trace Formula which provides the character of the representation
in terms of orbital integrals  
over conjugacy classes of the ad\`elic group. 
Furthermore, we have focused on the anisotropic conjugacy classes and expressed their orbital integrals in terms of $\kappa$-twisted orbital integrals
over stable conjugacy classes in $p$-adic groups. 

It is not too much of a stretch to argue that the $\kappa$-stable orbital integrals  $\CO^\kappa_{\gamma_0}(\varphi)$ are more appealing than the basic orbital integrals  $\CO_{\gamma}(\varphi)$ since 
 their dependence on the conjugacy
 classes within a stable conjugacy class is through a character rather than a specific choice of conjugacy class.
 This is an early manifestation of the motivic, or universal algebraic, nature of the $\kappa$-orbital integrals.
 But of course, aesthetics aside, Fourier inversion tells us we can go back and forth between the two, and so in some sense
 we have not accomplished very much.

Thus perhaps we have made a Faustian bargain: we have traded the evident geometric structure of basic orbital integrals  $\CO_{\gamma}(\varphi)$
for the representation theoretic structure of $\kappa$-orbital integrals $\CO^\kappa_{\gamma_0}(\varphi)$.
With our original aim  to compare trace formulas for different groups, one could even worry that we have made things more difficult rather than less. Indeed, one could argue that what we have done ``is obviously useless, because the term $\CO^\kappa_{\gamma_0}(\varphi)$
is still defined in terms of $G$" rather than some other group~\cite{harrisintro}.

But now we have arrived in the neighborhood of the Fundamental Lemma. It is the lynchpin in Langlands's theory of endoscopy
which relates  $\kappa$-orbital integrals to stable orbital integrals on other groups.
The theory of endoscopy (for which we recommend the original papers of Kottwitz~\cite{ko84, ko86})
has many facets, but at its center is the following question:
\begin{quote}
{\em For a local field $F$, given an element $\gamma_0\in G(F)$, and a compatible character $\kappa\in T^\vee$, 
on what group $H$ should we try to express the $\kappa$-orbital integral $\CO^\kappa_{\gamma_0}(\varphi)$ 
as a stable orbital integral?}
\end{quote}
The answer is what ones calls the {\em endoscopic group} associated to the given data.
At first pass, it is a very strange beast, neither fish nor fowl. But the Fundamental Lemma is what confirms it is the correct notion.

There is a great distance between the intuitive idea of an endoscopic group and the minimal notions one
needs to at least spell out the Fundamental Lemma. Most of the technical complications devolve from the intricacy of Galois descent
for  quasi-split groups. So it seems useful, though less efficient, to first explain the basic notions
assuming all groups are split (Definition~\ref{def split}), and then add in the necessary bells and whistles for quasi-split groups
(Definition~\ref{def quasisplit}).


\subsubsection{Split groups}
We begin with a reminder  of the ``combinatorial skeleton"
of a split reductive group given by its root datum. We will always equip all split reductive groups $G$ with a {\em pinning} consisting of a Borel subgroup $B\subset G$, split maximal torus $T\subset B$, and basis vectors for the resulting simple positive root spaces.
This has the effect of providing a canonical splitting
$$
\xymatrix{
1 \ar[r] & \Inn(G) = G/Z(G) \ar[r] & \Aut(G)\ar[r] & \Out(G) \ar[r]  \ar@/_/[l] & 1
}
$$
since the automorphisms of $G$ preserving the pinning map isomorphically to the outer automorphisms of $G$.

\begin{example} Take $G= SL(n)$ and $T \subset SL(n)$ the split maximal torus of diagonal matrices of determinant one. 

Then the symmetric group $\Sigma_n$ acts simply transitively on the  possible Borel subgroups $B\subset G$ satisfying $T\subset B$. Let us choose $B\subset SL(n)$ to consist of upper-triangular matrices of determinant one. 

The resulting simple positive root spaces can be identified with the $n-1$ super-diagonal matrix entries (directly above the diagonal). Let us choose the basis given by taking the element $1$ in each simple positive root space.

The outer automorphisms $\Out(SL(2))$ are trivial, but when $n>2$, the outer automorphisms $\Out(SL(n))$ are the group $\BZ/2\BZ$. The above pinning realizes the nontrivial outer automorphism
as the automorphism given by 
$$
\xymatrix{
A \ar@{|->}[r] & M (A^{-1})^\tau M^{-1},
&
A\in SL(n),
}
$$
where $\tau$ denotes the transpose operation, and  $M$ is the antidiagonal matrix with $M_{i, n-i+1} = 1$ when $i< n/2$, and $M_{i, n-i+1} = -1$ when $i\geq  n/2$.

\end{example}

\begin{defn}
(1) A (reduced) {\em root datum} is an ordered quadruple
$$
\Psi = (X, \Phi, X^\vee, \Phi^\vee)
$$ 
of the following data:
\begin{enumerate}
\item $X, X^\vee$ are finite rank free $\BZ$-modules in duality by a pairing 
$$
\xymatrix{
\langle , \rangle: X\times X^\vee\ar[r] & \BZ
}
$$
\item $\Phi, \Phi^\vee$ are finite subsets of $X, X^\vee$ respectively in fixed bijection
$$
\xymatrix{
\alpha \ar@{<->}[r] & \alpha^\vee
}$$
\end{enumerate}
We will always assume that our root data are reduced in the sense that if $\alpha \in \Phi$, then $c\alpha\in \Phi$ if and only if $c=\pm 1$.

The data must satisfy the following properties:
\begin{enumerate}
\item[(a)] $\langle \alpha, \alpha^\vee\rangle = 2$,
\item[(b)] $s_\alpha(\Phi) \subset \Phi, s_\alpha^\vee(\Phi^\vee)\subset \Phi^\vee$, where
$$
\begin{array}{ll}
s_\alpha(x) = x - \langle x, \alpha^\vee\rangle \alpha, 
&
x\in X, \alpha \in \Phi,\\
s^\vee_\alpha(y) = y - \langle \alpha, y\rangle \alpha^\vee, 
&
y\in X^\vee, \alpha \in \Phi.
\end{array}
$$
\end{enumerate}

The {\em Weyl group} $W_\Psi$ of the root datum is the finite subgroup of $GL(X)$  generated by 
the reflections $s_\alpha$, for $\alpha\in \Phi$. 

(2) A  {\em based root datum} is an ordered sextuple
$$
\Psi = (X, \Phi, \Delta, X^\vee, \Phi^\vee, \Delta^\vee)
$$ 
consisting of a root datum $ (X, \Phi, X^\vee, \Phi^\vee)
$
together with a choice of subsets
$$
\Delta\subset \Phi, \Delta^\vee \subset \Phi^\vee
$$
satisfying the following properties:
\begin{enumerate}
\item[(a)] the bijection $\Phi\longleftrightarrow \Phi^\vee$ restricts to a bijection $\Delta\longleftrightarrow \Delta^\vee$,
\item[(b)] there exists an element $v\in X$ with trivial stabilizer in $W_\Psi$, for which we have
$$
\Delta^\vee = \{ \alpha^\vee\in \Phi^\vee | \langle v, \alpha^\vee\rangle > 0\}
$$
\end{enumerate}

\end{defn}

To a split reductive group $G$ with a Borel subgroup $B\subset G$, and maximal torus $T\subset B$, one associates the based root datum
$$
\Psi(G) = (X_*, \Phi_G, \Delta_G, X^*, \Phi_G^\vee, \Delta^\vee_G)
$$
consisting of the following:
\begin{itemize}
\item $X_* = X_*(T) = \Hom(GL(1), T)$ the cocharacter lattice, 
\item $X^* = X^*(T) =  \Hom(T, GL(1))$ the character lattice, 
\item $\Phi_G\subset X_*$ the coroots, 
\item $\Phi_G^\vee\subset X^*$ the roots,
\item $\Delta_G\subset \Phi_G$ the simple coroots, and
\item $\Delta^\vee_G\subset \Phi_G^\vee$ the simple roots.
\end{itemize}

The Weyl group $W_{\Psi(G)}$ coincides with the usual Weyl group $W_G = N_T/T$.

Here is a key motivation for the notion of based root data.

\begin{thm} Fix a field $k$.

(1) Every based root datum $\Psi$ is isomorphic to the based root datum $\Psi(G)$ of some split reductive group $G$,
defined over $k$, and equipped
with a pinning.

(2) The automorphisms of the based root datum $\Psi(G)$ are isomorphic to the outer automorphisms of $G$,
or equivalently, the automorphisms of $G$, as an algebraic group defined over $k$, that preserve its pinning.
\end{thm}

The combinatorial classification of groups finds ubiquitous application, and is further justified by many natural occurrences
of related structures such as Dynkin diagrams. But one of its initially naive but ultimately deep implications is the evident
duality for reductive groups coming from the duality of root data.
It generalizes the very concrete duality for tori we have seen earlier which assigns to a split torus $T = \Spec k[X^*(T)]$ the dual complex torus 
$T^\vee = \Spec \BC[X_*(T)]$.

\begin{defn} Let $G$ be a split reductive group
with based root datum
$$
\Psi(G) = (X_*, \Phi_G, X^*, \Phi_G^\vee, \Delta_G, \Delta^\vee_G).
$$
The Langlands dual group $G^\vee$ is the split reductive complex algebraic group with
dual  based root datum 
$$
\Psi(G^\vee) = (X^*, \Phi_G^\vee,  \Delta^\vee_G, X_*, \Phi_G,  \Delta_G).
$$
\end{defn}

\begin{remark}
We have stated the duality asymmetrically, where $G$ is defined over some field $k$, but the dual group $G^\vee$ is always a complex algebraic group. Observe that such asymmetry arose for tori when we described complex characters. In general, it stems from the fact that our automorphic representations are complex vector spaces.

For a group $G$ with Langlands dual group $G^\vee$, the maximal torus  $T^\vee\subset G^\vee$ is the dual of the maximal torus  $T\subset G$, the Weyl groups $W_G$ and $W_{G^\vee}$ coincide, the outer automorphisms $\Out(G)$ and $\Out(G^\vee)$ coincide, and the roots of $G$
are the coroots of $G^\vee$ and vice-versa. If $G$ is simple, then so is $G^\vee$, and if in addition $G$ is complex, $Z(G)\simeq \pi_1(G^\vee)$
and vice-versa (generalizations of the last assertion are possible but one has to be careful  to compare potentially different kinds of groups).
\end{remark}

\begin{example}
The following are pairs of Langlands dual groups:
$GL(n) \longleftrightarrow GL(n)$,
$SL(n) \longleftrightarrow PGL(n)$,
$SO(2n+1) \longleftrightarrow Sp(2n)$,
 $SO(2n) \longleftrightarrow SO(2n)$.
\end{example}

Although the above definition is concrete, there is a deep mystery in passing from a group to root data,
dual root data, 
and then back to a group again. Commutative and combinatorial structures are the only things
which can easily cross the divide.\footnote{ There are many hints in quantum field theory of ``missing" higher-dimensional objects
which can be specialized on the one hand to reductive groups, and on the other hand to their root data. But until
they or their mathematical analogues are understood in some form, the relation of group to root data and hence to dual group will likely remain mysterious.}

Now we arrive at the notion of endoscopic group in the context of split groups.

\begin{defn}
Let $G$ be a split reductive algebraic group with split maximal torus $T\subset G$.

(1) {\em Split endoscopic data} is an element $\kappa \in T^\vee \subset G^\vee$.

(2) Given  split endoscopic data $\kappa\in T^\vee$, the associated {\em split endoscopic group} of $G$ is the split reductive algebraic group $H$ whose Langlands dual group $H^\vee$
is the connected component of the centralizer $G^\vee_\kappa \subset G^\vee$ of the element $\kappa$.
\end{defn}

It follows immediately that $T$ is also a maximal torus of $H$ and the coroots $\Phi_H$  are a subset of the coroots $\Phi_G$.
More precisely, the element $\kappa\in T^\vee = \Hom(X_*(T), GL(1))$ can be evaluated on $X_*(T)$, and in particular on $\Phi_G \subset X_*(T)$,
and the coroots $\Phi_H$ are given by the kernel
$$
\Phi_H = \{ \alpha\in \Phi_G | \kappa(\alpha) = 1\}.
$$
This immediately implies
that the roots $\Phi^\vee_H$ are the corresponding subset of the roots $\Phi_G$,
and  the Weyl group $W_H$ is a subgroup of the Weyl group $W_G$. 
{\em But  this by no means implies that $H$
is anything close to a subgroup of $G$}.

\begin{example} We will work with the split groups defined in Example~\ref{classical groups}.

Take $G=Sp(2n)$ the symplectic group so that $G^\vee = SO(2n+1)$ the odd orthogonal group.
Recall that the diagonal matrices inside of  $SO(2n+1)$ furnish a split maximal torus $T^\vee$. Take the element
$$
\kappa=
\left[
\begin{array}{ccc}
1 & 0  \\
 0 & -I_n
\end{array}
\right] \in T^\vee
$$
with centralizer $O(2n)$ which is disconnected with connected component  $H^\vee=SO(2n)$.
Taking the Langlands dual of $H^\vee$ gives the endoscopic group $H = SO(2n)$.

One can check that there is  no nontrivial homomorphism from $SO(2n)$ to $Sp(2n)$
when $n\geq 3$. (When $n=2$, there are homomorphisms, but they do not induce the correct maps on root data.)
\end{example}

\subsubsection{General story: quasi-split groups}
Now we will pass to the general setting of quasi-split groups. The inexperienced reader could skip this material
and still find plenty of interesting instances of the Fundamental Lemma to ponder.
Throughout what follows,  let $F$ be a non-Archimedean local field with residue field $f$. 

\begin{defn} 

(1)A field extension $E/F$ with residue field $e/f$ is said to be {\em unramified}
if the degrees satisfy $[E:F] = [e:f]$. 

(2) A reductive algebraic group $G$ defined over a local field $F$ is said to be {\em unramified} if the following hold:
\begin{enumerate}
\item[(a)] $G$ is quasi-split,
\item[(b)]  $G_{F^{un}}$ is split for the maximal unramified extension $F^{un}/F$.
\end{enumerate}

\end{defn}

Unramified groups $G$ are combinatorial objects, classified by the split reductive group $G_{F^{un}}$
together with
%
the Galois descent homomorphism
$$
\xymatrix{
\rho_G:Gal(F^{un}/F)\ar[r] &  \Out(G)
}
$$
which is completely determined by its value on the Frobenius automorphism.

%
%
%
%

To state the general notion of unramified endoscopic group, we will need the following standard constructions.
Given a split reductive group $G$, and an element $\kappa\in T^\vee \subset G^\vee$, recall that we write $H^\vee$ for the connected component
of the centralizer $G^\vee_\kappa\subset G^\vee$ of the element $\kappa$.
We will write $\pi_0(\kappa)$ for
the component group of the centralizer of the element $(\kappa, e)$ of the semi-direct product $G^\vee \rtimes \Out(G^\vee)$.
There are canonical maps
$$
\xymatrix{
 \Out(H^\vee) & \ar[l]_-{\pi_{H^\vee}} \pi_0(\kappa) \ar[r]^-{\pi_{G^\vee}} & \Out(G^\vee)
}
$$
$$\xymatrix{
W_{H^\vee}\rtimes \Out(H^\vee) \ar[r] & W_{G^\vee}\rtimes \Out(G^\vee) 
}
$$
where the latter is compatible with the actions on $T^\vee$ and projections to $\pi_0(\kappa)$, 
and extends the canonical inclusion $W_{H^\vee} \to W_{G^\vee}$.

\begin{defn} Let $G$ be an unramified reductive group, so in particular quasi-split with Borel subgroup $B\subset G$,
and (not necessarily split) maximal torus $T\subset B$.

(1)  Unramified {\em endoscopic data} is a pair $(\kappa, \rho_\kappa)$ consisting of an element $\kappa\in T^\vee \subset G^\vee$, and a homomorphism 
$$
\xymatrix{
\rho_\kappa: Gal(F^{un}/F) \ar[r] & \pi_0(\kappa)
& \mbox{ such that }  \rho_G = \pi_{G^\vee}\circ \rho_\kappa.
}
$$

(2) Given unramified {endoscopic data} $(\kappa, \rho_\kappa)$, the associated unramified {\em endoscopic group} of $G$ is the unramified reductive group $H$ defined over $F$ constructed as follows. Recall that the split endoscopic group  
$H^{spl}$ associated to the element $\kappa$ is the connected component of the centralizer $G_\kappa^\vee\subset G^\vee$. 
The endoscopic group $H$ is the form of $H^{spl}$ defined over $F$ by the
the  Galois descent homomorphism
$$
\xymatrix{
\rho_{H} = \pi_{H^\vee}\circ \rho_\kappa: Gal(F^{un}/F) \ar[r] & \Out(H^{spl}).
}
$$
\end{defn}

\begin{example}[Endoscopic groups for $SL(2)$]
There are three unramified endoscopic groups for $SL(2)$. Recall that the dual group is $PGL(2)$.
\begin{enumerate}
\item $SL(2)$ itself with $\kappa = e\in T^\vee$ the identity, and $\rho_\kappa$ trivial.
\item $GL(1)$ with $\kappa \not = e \in T^\vee$ any nontrivial element, and $\rho_\kappa$ trivial.
\item $U(1, E/F)$ with $\kappa = -e \in T^\vee$ the square-root of the identity, and $\rho_\kappa$ the non-trivial
map to $\pi_0(\kappa) \simeq \BZ/2\BZ$. 
\end{enumerate}
The first two are split, but the last is not.
\end{example}



\subsubsection{Statement of Fundamental Lemma}

Finally, we arrive at our destination.

The Fundamental Lemma relates  the $\kappa$-orbital integral
 over a stable conjugacy class in the group $G$
with the stable orbital integral  over a stable conjugacy class in 
an endoscopic group $H$.
Such an idea should lead to some immediate confusion: {\em the orbital integrals to be compared are distributions on different groups, so  to compare them
we must also have some correspondence of test functions.}

There is a deep and intricate theory of transferring test functions of which the Fundamental Lemma is in some sense the simplest
and thus  most important instance.
It states that in the most hospitable situation, the most simple-minded transfer of the simplest test functions leads to a good comparison of orbital integrals. There are many variations (twisted, weighted, ...) of the Fundamental Lemma,
but  the most important are now understood thanks to reductions to the Fundamental Lemma or extensions of the ideas of its proof.

\medskip

Fix an unramified group $G$ defined over $F$.
The fact that $G$ is unramified implies that it is the localization of a smooth affine group scheme $\CG$  defined over the ring of integers $\CO_F$ whose special fiber over the residue field $k$ is connected reductive.

\begin{defn}
%

A maximal compact subgroup $K\subset G(F)$ is said to be {\em hyperspecial} if there is a smooth affine group
scheme $\CG$ defined over  $\CO_F$ such that 
\begin{enumerate}
\item  $\CG_F = G$, 
\item $\CG(\CO_F) = K$, and 
\item $\CG_k$ is connected reductive.
\end{enumerate}
\end{defn}

\begin{lemma} A reductive algebraic group $G$ defined over a local field $F$ is unramified if and only if $G(F)$ contains
a hyperspecial maximal compact subgroup.
\end{lemma}

Now an endoscopic group $H$ is not a subgroup of $G$.
Rather we must content ourselves with
 the relationship of characteristic polynomials
$$
\xymatrix{
H \ar[d]_-{\chi_H} & G\ar[d]^-{\chi_G}\\
T/W_H \ar[r]^-\nu & T/W_G
}
$$

This provides a relationship between stable conjugacy classes as follows. Given a parameter $a_H\in (T/W_H)(F)$, we can consider its transfer $a_G = \nu(a_H) \in (T/W_G)(F)$.
Even if $a_H\in (T/W_H)(F)$ is regular, since the inclusion $W_H \subset W_G$ is not an isomorphism (except when $H= G$),
the transfer $a_G\in (T/W_G)(F)$ might not be regular.

\begin{defn}
 A parameter $a_H\in (T/W_H)(F)$ is said to be {\em $G$-regular} if it and its transfer $a_G = \nu(a_H) \in (T/W_G)(F)$ are both regular.
\end{defn}

Now given a  {\em $G$-regular} parameter $a_H\in (T/W_H)(F)$ with transfer $a_G =\nu(a_H)\in (T/W_G)(F)$, we have the stable conjugacy class 
$$[\gamma_H]_{st} = \chi_H^{-1}(a_H) \subset H
 \quad\text{ and its transfer }\quad
[\gamma_G]_{st} = \chi_G^{-1}(a_G) \subset G.
 $$

\begin{fundlemma}[\cite{ngoFL}]
Let $F$ be a local field. 

Let $G$ be an unramified group defined over $F$.

Let $H$ be an unramified endoscopic group for $G$ associated  to endoscopic data $(\kappa, \rho_\kappa)$.

Let $K_G\subset G(F), K_H \subset H(F)$ be hyperspecial maximal compact subgroups.

Then 
for  $G$-regular $a_H \in T/W_H$ and transfer $a_G = \nu(a_H)\in (T/W_G)(F)$,
we have an equality
$$
\CS\CO_{\gamma_H}(1_{K_H}) =   \Delta(\gamma_H, \gamma_G) \CO^\kappa_{\gamma_G}(1_{K_G})
$$
where $\gamma_H \in \chi^{-1}_H(a_H)$, $\gamma_G \in \chi^{-1}_G(a_G)$, and  $\Delta(\gamma_H, \gamma_G)$ is the transfer factor (which shall not be defined here).
\end{fundlemma}

\begin{remark}
A precise formulation of transfer factors first appears in Langlands's joint work with 
 D. Shelstad~\cite{langlandsshelstad}.
 The transfer factor $\Delta(\gamma_H, \gamma_G) $ accounts for the ambiguity that the $\kappa$-orbital integral 
$\CO^\kappa_{\gamma_G}(\varphi)$ depends on the choice of lift $\gamma_G \in \chi_G^{-1}(a_G)$. By definition, the stable orbital integral $\CS\CO_{\gamma_H}(\varphi) $ is an invariant of $a_H = \chi_H(\gamma_H)$.
It is worth mentioning that if $G$ is split and the derived group $G_{der} =[G,G]$ is simply-connected, then the Steinberg section provides a distinguished lift $\gamma_G =  \sigma(a_G).$
\end{remark}

\begin{remark}
There is an analogous ``Fundamental Lemma" for Archimedean local fields,
resolved long ago by D. Shelstad~\cite{shelstad}, which one also needs for applications of the trace formula. 
The example at the beginning of the Introduction fits into this  Archimedean part of the theory.
\end{remark}
%

It should be apparent that one can formulate a Lie algebra variant of the Fundamental Lemma. Namely, let $\fg$ be the Lie algebra of $G$, and $\fh$  the Lie algebra of an endoscopic group. Then one can
 replace the stable conjugacy classes of the group elements  $\gamma_H \in H(F)$ and $\gamma_G \in G(F)$ with those of 
 Lie algebra elements $\xi_H \in \fh(F)$ and $\xi_G \in \fg(F)$. 
 
 In fact, the situation simplifies further in that each stable conjugacy class now has a canonical element.
To see this, observe that stable conjugacy classes in the Lie algebra $\fg$ are the fibers of the Chevalley morphism 
$$
\xymatrix{
\chi:\fg \ar[r] & \ft/W   = \Spec  k[\ft]^W. 
}
$$
For the choice of a pinning, there is the Kostant section
$$
\xymatrix{
\sigma: \ft/W  \ar[r] & \fg.
}
$$
This is completely general, 
unlike for the group $G$, where the Steinberg section could be defined only when the derived group $G_{der} = [G, G]$
is simply-connected.

Thus with the assumptions of the Fundamental Lemma stated above, 
the Lie algebra variant takes the form of an identity
$$
\CS\CO_{a_H}(1_{\fh(\CO_F)}) =  \CO^\kappa_{a_G}(1_{\fg(\CO_F)})
$$
where we index the stable orbital integral by the parameter $a_H\in (\ft/W_H)(F)$, and we base the $\kappa$-orbital integral at the image of the Kostant section applied to the transfer $a_G = \nu(a_H)\in (\ft/W_G)(F)$.
In particular, the distinguished base point  obviates the need for any transfer factor.

An important theorem of  Waldspurger~\cite{wallie} asserts that the Lie algebra variant of the Fundamental Lemma  implies
the original statement. And it is the Lie algebra variant  which Ng\^o proves, and which we will
turn to in the next section.



\section{Geometric interpretation of Fundamental Lemma}\label{sectlocal}

%
%

In retrospect, the search for a proof of the Fundamental Lemma turns out to be the search for a setting where the powerful
tools of algebraic geometry -- Hodge theory, Lefschetz techniques, sheaf theory, homological algebra -- could be brought to bear on the problem.
As we have seen, the Fundamental Lemma is an analytic assertion about integrals of characteristic functions
on $p$-adic groups.
But these functions turn out to be of motivic origin, and hence amenable to the deep mid-to-late 20th century
synthesis of algebraic geometry and algebraic topology.

Here is a historical antecedent worth keeping in mind (which is deeply intertwined with the Fundamental Lemma and its proof). 
The
classical Riemann hypothesis 
that all non-trivial zeros of the Riemann zeta function have real part $1/2$
is a difficult problem. Even more dauntingly, it admits a well-known generalization to any number field $F$ for which 
the classical version is the case of the rational numbers $\BQ$.  Nevertheless, an analogous Riemann hypothesis for function fields of curves, formulated by Artin, was proved for elliptic curves by Hasse, and then for all genus curves by Weil. 
It involves counting the number of points of the curve over finite fields. The basic case of the projective line
is completely trivial. 
The Weil conjectures are a vast generalization to all algebraic varieties. They were a prominent focus
of Grothendieck and Serre, and eventually established by Deligne. At the heart of the proof is the interpretation
of  counting points in terms of cohomology with Galois actions. So in the end, the point counts have much more structure
than one might have thought. 

\subsection{Motivic origin}

Though the Fundamental Lemma is a local statement, its motivation comes from global questions about a number field.
Without Langlands functoriality and the Arthur-Selberg Trace Formula in mind, it would be hard to arrive at the Fundamental Lemma as a reasonable assertion. 
Nevertheless, once we dispense with motivation, the concrete problem to be solved involves only local fields.

Remarkably,  the analogy between number fields and function fields of curves leads to precise mathematical comparisons
between their completions. It is hard to imagine that all of the intricate arithmetic of a number field could be found
 in the geometry of an algebraic curve. For example,
the structure of the rational numbers $\BQ$ is far more complicated than that of 
 the projective line $\BP^1$.   But it turns out that important structures of the $p$-adic fields $\BQ_p$ can be found
in  the Laurent series fields $\BF_p((t))$. Naively, though $\BQ_p$ is of characteristic zero and $\BF_p((t))$ is of characteristic
$p$, they share the formal structure that each is a complete local field with residue field $\BF_p$.
In particular, one can transport the statement of the Fundamental Lemma from the setting of $p$-adic fields to
 Laurent series fields.

The following liberating theorem of Waldspurger opens the door to geometric techniques.

\begin{thm}[\cite{walchar}]
The Lie algebra variant of the Fundamental Lemma in the equal characteristic or geometric case of $\BF_p((t))$ implies the Lie algebra variant of the Fundamental Lemma
in the unequal characteristic or arithmetic case of $\BQ_p$.
\end{thm}

\begin{remark}
Waldspurger proves the above assertion for sufficiently large residual characteristic $p$. Thanks to previous work of Hales~\cite{hales}, this suffices to establish the Fundamental Lemma for arbitrary residual characteristic.
\end{remark}

Waldspurger's proof is a tour de force of representation theory. 
It provides a detailed analysis of constructions which are natural but specific
 to reductive groups. 
 
Cluckers, Hales, and Loeser~\cite{chl} have discovered a completely independent proof of the above result.
The context of their arguments is mathematical logic, specifically model theory. While this is unfamiliar
to many, it is very appealing in that it gets to the heart, or motivic truth, of the independence of characteristic.
Very roughly speaking, the arguments take seriously the idea that integrals contain more information than their numerical values.
Rather, they should be thought of as universal linear expressions for their cycles of integration weighted by their
integrands. 
It applies generalizations of Cluckers-Loeser~\cite{cl1, cl2}, building on work of Denef-Loeser~\cite{dl},  of the classical 
Ax-Kochen-Ersov Theorem that given  $\varphi$ a first order sentence (a formula with no free variables) in the language of rings, for almost all prime numbers $p$, the sentence $\varphi$ is true in $\BQ_p$ if and only if it is true in $\BF_p((t))$.



\subsection{Affine Springer fibers}

In the geometric setting of Laurent series fields, thanks to the Weil conjectures, 
the Fundamental Lemma takes on a topological form.
It comes down to comparing the cohomology of affine Springer fibers for endoscopic groups.
We will explain what affine Springer fibers are and some beautiful techniques of broad applicability
developed to understand their cohomology.

Affine Springer fibers were introduced by Kazhdan-Lusztig~\cite{KL} as natural generalizations of  Grothendieck-Springer fibers.
The latter play a fundamental role in Springer's theory of Weyl group representations, as well as Lusztig's theory
of character sheaves. If we take the viewpoint that Grothendieck-Springer fibers are essential to
the characters of groups over a finite field, then it is not surprising that affine Springer fibers figure prominently in
the characters of $p$-adic groups.

For simplicity, and to appeal to topological methods, we will work here with base field the complex numbers $\BC$.
We write $\CK = \BC((t))$ for the Laurent series field, and  $\CO = \BC[[t]]$ for its ring of integers.

\begin{remark}
Since $\BC$ is algebraically closed, all unramified groups over $\CK$ will be split, and so we will not be discussing
the most general form of the Fundamental Lemma. We hope the broader intutions available in this setting will compensate for the sacrifice of generality.
\end{remark}

Let $G$ be a complex reductive group with Lie algebra $\fg$.
Consider the so-called affine or loop group
$G_\CK = G(\BC((t)))$, and its subgroup of arcs $G_\CO= G(\BC[[t]])$. 
The quotient $Gr_G = G_\CK/G_\CO$ is called the affine Grassmannian of $G$. It is an increasing union of 
projective varieties indexed by the natural numbers. The group $G_\CK$ naturally acts on $Gr_G$ by multiplication,
and we can think of its Lie algebra $\fg_\CK$ as acting infinitesimally by vector fields.

\begin{defn}
For an element $\xi\in \fg_\CK$ of  the Lie algebra of the loop group $G_\CK$, the affine Springer fiber $Gr_G^\xi$ is the fixed-points of the vector field $\xi$ acting on the affine Grassmannian
$Gr_G$. It admits the concrete description
$$
Gr_G^\xi = \{ g \in G_\CK | \on{Ad}_g (\xi)  \in \fg_\CO \}/G_\CO
$$
where $\fg_\CO$ denotes the Lie algebra of $G_\CO$.
\end{defn}

Although affine Springer fibers are not finite-type schemes, when the parameter $\xi$ is regular semisimple, their underlying reduced schemes are reasonable. We will see this by examining the natural symmetry groups acting upon them.

First, observe that the Springer fiber only depends upon the $G_\CK$-conjugacy class of $\xi$, and if $\xi$ is not conjugate to any $\xi' \in \fg_{\CO}$ then the Springer fiber will be empty. So there is no loss in assuming $\xi\in \fg_{\CO}$ from the start. Moreover, any regular semisimple $\xi\in \fg_\CO$ is $G_\CK$-conjugate to a unique $\xi_0 =\sigma(a) \in \fg_{\CO}$ in the image of the Kostant slice 
$$
\xymatrix{
\sigma:(\ft/W)_{\CO} \ar[r] & \fg_{\CO}.
}
$$
So in what follows, we will assume that $\xi =\sigma(a)$ for a fixed regular element $a\in (\ft/W)_\CO$.

Now, let $T = (G_{\CK})_\xi \subset G_\CK$ denote the (possibly ramified) torus given by the centralizer of $\xi$.
Let $S \subset T$ be its maximally unramified subtorus. In other words,  if we write $X_*(T_{\ol \CK})$ for the cocharacters of $T$,
then the cocharacters of $S$ are the Galois-invariants 
$$
X_*(S) = X_*(T_{\ol \CK})^{Gal(\ol \CK/\CK)}.
$$

A key observation is that $T$ canonically extends to a smooth commutative group-scheme $J$ defined over $\CO$.
Namely,  over the regular locus  $\fg^{reg} \subset \fg$, we have the smooth commutative group-scheme of centralizers $I \to \fg^{reg}$,
and we can form the fiber product
$$
J = \Spec\CO \times_\fg I
$$
over the base point $\xi= \epsilon(a):\Spec \CO \to \fg^{reg} \subset \fg$.
 
 We will write $\Lambda_\xi$
for the coweight lattice $X_*(S)$, and $\CP_\xi$ for the affine Grassmannian $J_{\CK}/J_{\CO}$.
The latter is an increasing union of 
not necessarily projective varieties. Since $J$ is commutative, $\CP_\xi$ is naturally a commutative group.
From the constructions, $\CP_\xi$ and hence also its subgroup $\Lambda_\xi$ naturally act on $Gr_G^\xi$.

\begin{prop}[\cite{KL}]
 Suppose $\xi\in\fg_\CK$ is regular semisimple. 
 
 Then the underlying reduced scheme $Gr_{G, red}^\xi$ of the affine Springer fiber is a countable union of projective
 irreducible components. Furthermore, the natural $\Lambda_\xi$-action on $Gr_{G, red}^\xi$ is free with quotient $Gr^\xi_{G, red}/\Lambda_\xi$ a projective variety.
\end{prop}

\begin{example} Take $G = SL(2)$, and consider the elements
$$
\xi_1 = 
\left[
\begin{array}{cc}
 0 & t \\
 1 & 0
\end{array}
\right],
\xi_2 = 
\left[
\begin{array}{cc}
 0 & t^2 \\
 1 & 0
\end{array}
\right],
\xi'_2 = 
\left[
\begin{array}{cc}
 t & 0 \\
 0 & -t
\end{array}
\right]
\in \fg_\CK.
$$
Then (at the level of reduced schemes) $Gr_G^{\xi_1}$ is a point, and $Gr_G^{\xi_2}$ is an infinite string of $\BP^1$'s attached one after the other at nodes. The lattice $\Lambda_{\xi_2}$ is rank one, and the quotient $Gr_G^{\xi_2}/\Lambda_{\xi_2}$ is a nodal elliptic curve.

To see things for $\xi_2$, it helps to note that $\xi_2$ and $\xi'_2$ are $G_\CK$-conjugate, and so one can calculate with $\xi_2'$.
On the other hand, $\xi_1$ is anisotropic and can not be diagonalized.

\end{example}

Now given a character 
$$
\kappa:\CP_\xi/\Lambda_\xi \to \BC^\times,
$$ 
we can consider the summand 
$$
H^*(Gr^\xi/\Lambda_\xi)_\kappa\subset 
H^*(Gr^\xi/\Lambda_\xi)
$$ 
consisting of  cocycles that transform by
the character $\kappa$
under the action of $\CP_\xi/\Lambda_\xi$.

Thanks to the following application of the Weil conjectures, the Fundamental Lemma involves identifying this cohomology
with the invariant cohomology of an affine Springer fiber  for an endoscopic group.

\begin{thm}
 Suppose $\xi\in\fg_\CK$ is regular semisimple. Then the $\kappa$-orbital integral $\CO^\kappa_{\xi}(1_{\fg_{\CO}})$
 can be recovered from the $\kappa$-summand $H^*(Gr^\xi/\Lambda_\xi)_\kappa$, or more precisely from its
 underlying motive.
\end{thm}

\begin{remark}
By underlying motive, we mean that one should properly work over a finite field with $\ell$-adic cohomology,
and use the Grothendieck-Lefschetz fixed point formalism to recover orbital integrals from traces of Frobenius.
\end{remark}

Unfortunately, there is no evident geometric relationship between the affine Springer fibers 
for a group and an endoscopic group. And any attempt to explicitly calculate the cohomology must face the fact that already for the symplectic
group $Sp(6)$, Bernstein and Kazhdan~\cite[Appendix]{KL} found an affine Springer fiber whose cohomology contains the motive 
of a hyperelliptic curve. 

%

  


\subsection{Equivariant localization}

Goresky-Kottwitz-MacPherson proposed the intriguing idea that the combinatorics of endoscopic groups
hint at a possible geometric mechanism for relating affine Springer fibers.  Recall that the roots of an endoscopic group are an explicit subset of the roots of the original group. For an affine Springer fiber with large toric symmetry, 
Goresky-Kottwitz-MacPherson recognized that its toric one-skeleton of zero and one-dimensional toric orbits was in fact completely encoded by the roots.
In turn, they also discovered that for a very general
class of varieties with large toric symmetries, their cohomology could be read off from their toric one-skeleta.

Let $X$ be a (possibly singular) projective variety equipped with an action of the torus $T = (\BC^\times)^n$.
For simplicity, we will fix an embedding $X\subset \BC\BP^N$, and assume the $T$-action is induced
by a linear $T$-action on $\BC^{N+1}$. Then there is a moment map 
$$
\xymatrix{
\mu:X\ar@{^(->}[r] &  \BC\BP^N \ar[r] & \ft^\vee
}
$$
which induces the corresponding infinitesimal action of $\ft$.
For a vector $v\in \ft$, regarded as a vector field on $\BC\BP^N$, the one-form $\omega(v, -)$
obtained by contracting with the Kahler form is given by the pairing
$
\langle d\mu(-), v\rangle 
$
Furthermore, we will assume the images of the fixed points $X^T \subset X$ are all distinct,
in particular, there are only finitely many.

\begin{defn}
Suppose $T = (\BC^\times)^n$ acts on $X$ with finitely many one-dimensional orbits. Let $X_1\subset X$ be the one skeleton
of fixed points and one-dimensional orbits. The moment graph $\Gamma_T(X)$ is defined to be the quotient
$$
\Gamma_T(X) = X_1/ T_c,
$$
where $T_c = (S^1)^n$ is the compact torus inside $T$.
The moment map descends to a canonical immersion
$$
\xymatrix{
\mu :\Gamma_T(X) \ar[r] & \ft^\vee.
}
$$
\end{defn}

The structure of the moment graph is not only the abstract graph, but also its immersion into $\ft^\vee$.
The image of any one-dimensional orbit  $\CO\subset X_1$ will be a line segment $\ell \subset \ft^\vee$
whose orthogonal in $\ft$ is the Lie algebra of the stabilizer of $\CO$.

\begin{example}
Fix a maximal torus $T\subset G$. For $\xi\in \ft_\CO\subset \fg_\CK$, 
the centralizer $Z_{G_\CK}(\xi)$ contains the affine torus $T_\CK$,
and hence the full coweight lattice $\Lambda = \on{Hom}(\BC^\times, T)$.
If $\xi$ is regular semsimple, the maximal torus $T\subset G$ acts on $Gr_G^\xi/\Lambda$ with finitely many one-dimensional orbits. The resulting moment graph $\Gamma_{T}(Gr_G^\xi/\Lambda)$ is an invariant of the root data of $G$.
\end{example}

The following equivariant localization theorem is not difficult but gives a beautiful combinatorial picture of the cohomology of $X$.
Its validity depends on the technical assumption that the $T$-action on $X$ is equivariantly formal. We will not explain 
what this means, but only mention that it follows from more familiar conditions such as the vanishing of the odd degree cohomology of $X$, or if the mixed Hodge structure on the cohomology of $X$ is in fact pure.

\begin{thm}[\cite{gkm98}]
Suppose the $T$-action on  $X$ is equivariantly formal.
Then the cohomology of $X$ is an invariant of the moment graph $\Gamma_T(X)$.
 \end{thm}
 
 This general viewpoint on the cohomology of spaces with torus actions has been very fruitful.
 Foremost, there is the original application of Goresky-Kottwitz-MacPherson to the Fundamental Lemma.
 
 \begin{corollary}[\cite{gkm04}]
 For $\xi\in \ft_\CO\subset \fg_\CK$ regular semisimple, if the cohomology of $Gr_G^\xi/\Lambda$ is pure, then the Fundamental Lemma holds for $\xi$.
 \end{corollary}
 
\begin{remark} Note that the corollary assumes $\xi$ lies in the loop algebra of $\ft$. This is a very restrictive condition on a regular semisimple element. In particular, it implies that its centralizer is an {\em unramified} torus. This starting point is what provides sufficient symmetries to apply equivariant localization to the affine Springer fiber.
 \end{remark}
 
Outside of the question of what to do with affine Springer fibers with less symmetry, the issue of purity is a formidable obstacle
to further progress in the local setting. Some cases were directly settled by Goresky-Kottwitz-MacPherson~\cite{gkm06}, but for a uniform understanding,
new ideas appeared necessary.


\section{Hitchin fibration and the search for purity}\label{sectglobal}

In this section, we will continue to work over the complex numbers $\BC$ in order to appeal to topological methods of broad familiarity. We will adopt all of the notation of the previous section, so for example we write $\CK = \BC((t))$ for the Laurent series field, and $\CO = \BC[[t]]$ for its ring of integers. Given the reductive group $G$, we write $G_\CK$ for its loop group, and $G_\CO$ for the subgroup of arcs.

\subsection{Grothendieck-Springer resolution}
Affine Springer fibers for the loop group $G_\CK$ are a natural generalization of Springer fibers for the original group $G$. 

For a nilpotent element $\xi\in \CN\subset \fg$ of  the Lie algebra of $G$, the Springer fiber $\CB^\xi$ is the fixed-points of the vector field $\xi$ acting on the flag variety $\CB$ of all Borel subalgebras $\fb\subset \fg$. 
They naturally arise
as the fibers of the Springer resolution of the nilpotent cone
$$
\xymatrix{
\mu_\CN: \tilde \CN \simeq T^*(G/B) = \{(\fb, \xi) \in \CB\times \CN | \xi\in\fb\} \ar[r] & \CN &
 (\fb, \xi) \ar@{|->}[r] & \xi
}
$$
Beginning with Springer's construction of Weyl group representations in their cohomology,  
Springer fibers are now ubiquitous in representation theory. 

A vital observation is that the Springer resolution may be extended to the so-called Grothendieck-Springer resolution
$$
\xymatrix{
\mu_\fg: \tilde\fg = \{(\fb, \xi) \in \CB\times \fg | \xi\in\fb\} \ar[r] & \fg &
 (\fb, \xi) \ar@{|->}[r] & \xi
}
$$
While Springer fibers for nilpotent elements $\xi \in \CN$ are difficult and important, Springer fibers for regular semisimple
elements $\xi \in \fg^{reg, ss}$ are finite and on their own uninteresting. So what have we accomplished by introducing
the Grothendieck-Springer resolution?
{\em We can now reduce questions about interesting Springer fibers to questions about dull Springer fibers.}

What technique can we use to relate the cohomology of the fibers of a map such as the Grothendieck-Springer resolution?
 Sheaf theory. The cohomology of the fibers of a map are precisely the local invariants of the derived pushforward
 of the constant sheaf along the map. Global results about the pushforward will imply local results about the cohomology
 of the fibers. For example, much of Springer's theory of Weyl group representations is encoded in the following statement
 (see for example ~\cite{Ginz, HK}).

 \begin{thm}\label{springer theory}
 The restriction of the derived pushforward $R\mu_{\fg*} \BC_{\tilde \fg}$ of the constant sheaf on $\tilde\fg$ to the open locus of regular semisimple elements $\fg^{reg, ss}\subset \fg$ is a local system with monodromies given by the regular representation of the Weyl group of $\fg$. The entire pushforward is the canonical intersection cohomology extension
  of this local system.
 \end{thm}

Returning to the loop group $G_\CK$, the affine Springer fibers $Gr_G^\xi$  are also fibers of an analogous map. But here the target
 is the infinite-dimensional Lie algebra $\fg_\CK$ and the fibers $Gr_G^\xi$ are not projective varieties. 
The powerful topological methods of algebraic geometry
do not (yet) extend to such a setting. 
Some new idea is needed  to proceed.

\subsection{Compactified Jacobians}
Laumon introduced a beautiful way to begin to study the relation of affine Springer fibers $Gr_G^\xi$ for varying
parameters $\xi\in \fg_\CK$. As we vary the parameter $\xi\in \fg_\CK$, the behavior of the affine Springer fibers
$Gr_G^\xi$ is far wilder than we might expect. 

\begin{example}\label{local node to cusp}
Take $G=SL(2)$, and consider the family of elements 
$$
\xi(\varepsilon) = 
\left[
\begin{array}{cc}
 0 & \varepsilon t^2 + t^3 \\
 1 & 0
\end{array}
\right] \in \fg_\CK,
\text{ for $\varepsilon\in \BC$}.
$$
When $\epsilon = 0$, (at the level of reduced schemes) the affine Springer fiber $Gr_G^{\xi(\varepsilon)}$ is simply $\BP^1$, but when $\varepsilon \not = 0$, it is an infinite string of $\BP^1$'s attached one after the other at nodes with symmetry lattice $\Lambda_{\xi(\varepsilon)} \simeq \BZ$.
\end{example}

Laumon recognized that the quotients $Gr_G^\xi/\Lambda_\xi$ which
interest us most in fact form reasonable families in a topological sense.

\begin{prop}[\cite{la06}]
Suppose $G= GL(n)$, and 
let $\xi\in \fg_\CK$ be regular semisimple. Then there is a complex  projective curve $C_\xi$ of genus zero, with a single singular
point, such that the quotient $Gr_G^\xi/\Lambda_\xi$ is homeomorphic to the compactified Jacobian
$$
\ol{{Jac}}(C_\xi)= \{\text{degree $0$, rank $1$ torsion free sheaves on $C_\xi$}\}.
$$
\end{prop}

This is a powerful insight: such curves $C_\xi$ form nice finite-dimensional families,
and hence so do their compactificed Jacobians $\ol{{Jac}}(C_\xi)$.

\begin{example}(compare with Example~\ref{local node to cusp}).
Take $G=SL(2)$, and consider the family of elements 
$$
\xi(\varepsilon) = 
\left[
\begin{array}{cc}
 0 & \varepsilon t^2 + t^3 \\
 1 & 0
\end{array}
\right] \in \fg_\CK,
\text{ for $\varepsilon\in \BC$}.
$$
Then we can take $C_{\xi(\varepsilon)}$ to be the family of singular elliptic curves $y^2 =\varepsilon t^2 + t^3$.
When $\epsilon = 0$, the curve $C_{\xi(0)}$ is a cuspidal elliptic curve, and when $\varepsilon \not = 0$, it is a nodal elliptic curve.
The compactified Jacobian of any elliptic curve, smooth or singular, is isomorphic to the curve itself. 
\end{example}

By deforming such curves, Laumon provides a natural geometric setting for the equivariant localization techniques
of Goresky-Kottwitz-MacPherson. In particular, he was able to deduce the Fundamental Lemma for unitary groups
contingent upon the  purity of the affine Springer fibers involved ~\cite{la}. 
The main obstruction for further progress remained the 
elusive purity on which all conclusions were conditional.

\subsection{Hitchin fibration}
Ng\^o recognized that Laumon's approach to affine Springer fibers via compactified Jacobians is a natural piece
of the Hitchin fibration. Although it might appear complicated at first glance, 
the Hitchin fibration is nothing more than the natural generalization of the Chevalley morphism 
$$
\xymatrix{
\chi:\fg \ar[r] & \ft/W   = \Spec  k[\ft]^W. 
}
$$
to the setting of algebraic curves. Recall that $\chi$ descends to the quotient $\fg/G$, and is also equivariant for $GL(1)$ acting on $\fg$ by  linear scaling, and on $(\ft/W)$ by the resulting weighted scaling.

Fix a smooth projective curve $C$. Although the constructions to follow are very general,  for simplicity we will continue to work over the complex numbers $\BC$, so in particular $C$ is nothing more than a compact Riemann surface. Let us imagine the Lie algebra $\fg$ varying as a vector bundle along $C$. In order to preserve its natural structures, we should insist that the transition functions of the vector bundle take values in $G$ acting by the adjoint representation. In other words,
 the twisting of the vector bundle should be encoded by a principal $G$-bundle $\CP$ and the vector bundle should take the form of the associated bundle 
 $$
 \fg_\CP =  \CP \times_G \fg.
 $$ 
We will write $\Bun_G$ for the moduli of principal $G$-bundles over $C$

Now suppose we are given a line bundle $\CL$ over $C$.
The $GL(1)$-action of  linear scaling on $\fg$ together with the line bundle $\CL$ provides the option to twist
 further and form the tensor product 
 $$
\fg_{\CP, \CL} = \fg_\CP \otimes_{\CO_C} \CL.
 $$
 Similarly, the  $GL(1)$-action of weighted scaling on $\ft/W$ together with the line bundle $\CL$ provides an affine bundle 
 $$
 (\ft/W)_{\CL} =  \CL \times_{GL(1)} \ft/W.
 $$

\begin{defn} Fix a smooth projective curve $C$ equipped with a line bundle $\CL$.

The total space $
\CM_G
$
of the Hitchin fibation 
is the moduli
of pairs called Higgs bundles of a principal $G$-bundle  $\CP$ over $C$, and a section of the twisted adjoint bundle
$$
 \varphi \in \Gamma(C,  \fg_{\CP, \CL})
$$
called a Higgs field.

The base $\CA_G$ of the Hitchin fibration is the space of possible eigenvalues
$$
\CA_G = \Gamma(C,  (\ft/W)_{\CL}).
$$
The Hitchin fibration is the pointwise unordered eigenvalue map
$$
\xymatrix{
\chi:\CM_G \ar[r] & \CA_G.
}
$$

The {Hitchin fibers} are the inverse images $\CM_G^a= \chi^{-1}(a)$ for parameters $a\in \CA_G$.
\end{defn}

\begin{remark}
Hitchin's original construction~\cite{hitchin} focused on the case where $\CL$ is the canonical line bundle $\omega_C$ of one-forms on $C$.
This is most natural from the perspective that then $\CM_G$ is the cotangent bundle to the moduli $\Bun_G$ of principal $G$-bundles,
and the Hitchin fibration is a complete integrable system.
 The choice of line bundle $\CL$ provides useful technical freedom, since by choosing $\CL$
ample enough, we can eliminate any constraint imposed by the global geometry of $C$.

\end{remark}

\begin{remark}
Given a line bundle $\CL$ and principal $G$-bundle $\CP$,
 we can find finitely many points of $C$ so that the restrictions of the bundles to their complement  are trivializable. Thus any point of $\CM_G$ gives an element of $\fg(F)$,
well-defined up to dilation and conjugation,
 where $F$ is the function field of $C$.
 Likewise any point of $\CA_G$ 
 gives an element of $(\ft/W)(F)$,
well-defined up to dilation.

In this way, we can import abstract definitions for algebraic groups defined over a field $F$ to the setting of the Hitchin fibration.
For example, a point of $\CM_G$ is said to be {regular}, {semisimple}, {regular semisimple},
or {anisotropic} if the generic value of its Higgs field is so as an element of $\fg(F)$.
Likewise  a point of $\CA_G$ is said to be generically {regular} if its generic value is a regular value of $(\ft/W)(F)$.
\end{remark}

\begin{example}[Hitchin fibration for $GL(n)$]\label{hitch gln}
Recall  that for $\fg \fl (n) = \End(\BC^n)$, the Chevalley morphism is simply the characteristic polynomial 
$$
\xymatrix{
\chi:\End(\BC^n) \ar[r] & \bigoplus_{k =1}^{n} L_k &
\chi(A) = \det(t \Id - A)
}
$$
where $L_k \simeq \BC$ denotes the one-dimensional vector space generated by the $k$th elementary
symmetric polynomial.

The total space $
\CM_{GL(n)}
$
is the moduli
of pairs of a rank $n$ vector bundle  $\CV$ over $C$, and  a twisted endomorphism
$$
 \varphi \in \Gamma(C, End_{\CO_C}(\CV) \otimes_{\CO_C} \CL).
$$

The base $\CA_{GL(n)}$ is the space of possible eigenvalues
$$
\CA_{GL(n)} = \bigoplus_{k =1}^{n} \Gamma(C, L_k \otimes_{\CO_C} \CL^{\otimes k})).
$$
The Hitchin fibration is the pointwise unordered eigenvalue map
$$
\xymatrix{
\chi:\CM_{GL(n)}\ar[r] & \CA_{GL(n)}.
}
$$

A parameter $a\in \CA_{GL(n)}$ assigns to each point $c\in C$ a collection of $n$ unordered eigenvalues.
The {\em spectral curve}  $C_a$ is the total space of this varying family over $C$. More precisely, it is the solution to the equation on the total space of $\CL$ given by the characteristic polynomial corresponding to $a$.
When the spectral curve $C_a$ is reduced, the Hitchin fiber 
is isomorphic to the compactified Picard
$$
\CM_{GL(n)}^a\simeq
\ol{{Pic}}(C_a)
$$
of rank $1$, torsion free sheaves on $C_a$.

\end{example}

\begin{example}[Hitchin fibration for $SL(2)$]\label{hitch sl2}

For $\fs \fl (2) =\{ A\in \fg\fl(2) | \Tr(A) = 0\}$, the Chevalley morphism is simply the determinant map
$$
\xymatrix{
\chi:\fs \fl (2)  \ar[r] & \BC &
\chi(A) = \det(A)
}
$$

The total space $
\CM_{SL(2)}
$
is the moduli
of pairs of a rank $2$ vector bundle  $\CV$ over $C$ with trivialized determinant $\wedge^2 \CV \simeq\CO_C$, and  a twisted traceless endomorphism
$$
 \varphi \in \Gamma(C, End_{\CO_C}(\CV) \otimes_{\CO_C} \CL),
 \qquad \Tr(\varphi) = 0 \in \Gamma(C, \CL).
$$

The base $\CA_{SL(2)}$ is the space of possible determinants
$$
\CA_{SL(2)} = \Gamma(C, \CL^{\otimes 2}).
$$
The Hitchin fibration is the pointwise determinant map
$$
\xymatrix{
\chi:\CM_{SL(2)}\ar[r] & \CA_{SL(2)}.
}
$$

We can regard a parameter $a\in \CA_{SL(2)}$ as an element of $\CA_{GL(2)}$ and hence assign to it a spectral curve $C_a$.
This will be nothing more than the two-fold cover $c_a:C_a \to C$ given by the equation $t^2+a$ on the total space of $\CL$.
When $a$ is not identically  zero,  $C_a$ is reduced, and the Hitchin fiber 
is isomorphic to the moduli
$$
\CM_{SL(2)}^a\simeq
\{ \CF \in \ol{{Pic}}(C_a) | \det(c_{a*}\CF) \simeq \CO_C\}.
$$
\end{example}


The importance of the Hitchin fibration in gauge theory, low dimensional topology, and geometric representation theory
can not be underestimated. Note that the moduli $\Bun_G$ of principal $G$-bundles is a precise analogue of the primary locally symmetric space of automorphic representation theory. Namely, if  $F$ denotes the function field of the smooth projective complex curve $C$, 
and $\BA_F$ its ad\`eles with ring of integers $\CO_F$, then
the moduli of principal $G$-bundles is isomorphic to the double coset quotient
$$
\Bun_G \simeq G(F)\backslash G(\BA_F)/G(\CO_F).
$$

So in some sense we have come full circle. Starting from questions about number fields, we arrived at $p$-adic local fields and the Fundamental Lemma. Then we translated the questions to Laurent series local fields, and finally
we will appeal to the geometry of function fields.


Now instead of considering affine Springer fibers, 
Ng\^o proposes that we consider Hitchin fibers.
To spell out more precisely the relation between the two, we will focus on their symmetries.

Recall that 
to a regular semisimple Kostant element $\xi \in \fg_{\CO}$, we associate a 
commutative group-scheme  $J$ over $\CO$, and the affine Grassmannian $\CP_\xi = J_{\CK}/J_{\CO}$
naturally acts on the affine Springer fiber $Gr_G^\xi$. 

Similarly, there is a global version of this construction which associates
to a generically regular element $a\in \CA_G$, a commutative group scheme $J$ over the curve $C$, and the moduli $\CP_a$
of principal $J$-bundles is a commutative group-stack which naturally acts on the Hitchin fiber $\CM_G^a$.
Here generically regular means that the value of the section $a$ is regular except at possibly finitely many points of $C$.

The following result of Ng\^o is a direct analogue of the ad\`elic factorization appearing in Formula~\eqref{adelic factorization}.

\begin{thm}\label{geometric adelic factorization}
Suppose $a\in \CA_G$ is generically regular. 

Let $c_i \in C$, for $i\in I$, be the finitely many points where $a$ is not regular, and let $D_i = \Spec \CO_i$ be the formal disk around $c_i$.

Consider the Kostant elements $\xi_i = \sigma(a|_{D_i}) \in \fg_{\CO_i}$.
Then there is a canonical
map inducing a topological equivalence
$$
\xymatrix{
\prod_{i\in I} Gr_G^{\xi_i} / \CP_{\xi_i} \ar[r]^-\sim & \CM_G^a/ \CP_a.
}
$$

 \end{thm}

\begin{example}[Example~\ref{hitch gln} of $GL(n)$ continued]

When the spectral curve $C_a$ of the parameter $a\in \CA_{GL(n)}$ is reduced, the symmetry group $\CP_a$ is precisely the Picard $Pic(C_a)$ of line bundles on $C_a$. Under the identification of the Hitchin fiber
$$
\CM_{GL(n)}^a\simeq
\ol{{Pic}}(C_a)
$$
with the compactified Picard, the action of $\CP_a\simeq Pic(C_a)$
 is simply tensor product.
\end{example}

\begin{example}[Example~\ref{hitch sl2} of $SL(2)$ continued]

When  the parameter $a\in \CA_{GL(n)}$ is nonzero so that $C_a$ is reduced, the symmetry group $\CP_a$ is the Prym variety
given by the kernel of the norm map
$$
\xymatrix{
\CP_a \simeq \ker \{N: Pic(C_a) \ar[r] & Pic(C)
\}
}
$$
Under the identification of the Hitchin fiber explained above,
the action of $\CP_a$
 is simply tensor product.
\end{example}

From the theorem and careful choices of the  parameter $a\in \CA_G$
and  the line bundle $\CL$, one can realize the cohomology of affine Springer fibers 
completely in terms of the cohomology of Hitchin fibers. Now we are in a finite-dimensional setting where the tools of algebraic geometry apply. Without developing any further theory, Laumon and Ng\^o~\cite{LN} were able to establish the necessary
purity to deduce the Fundamental Lemma for unitary groups.


\subsection{Ng\^o's Support Theorem} 

From our preceding discussion, we conclude that we can replace the study of affine Springer fibers with that of Hitchin fibers.
Unlike the somewhat discontinuous behavior of affine Springer fibers, the Hitchin fibration is a highly structured family. To understand its cohomology,
Ng\^o introduces the following general notion of an abelian fibration. 
What he proves about them will no doubt find much 
further application. 
For simplicity, we will continue to work over the complex numbers $\BC$,
though the notions and results make sense in great generality.

%
%

\begin{defn} A {\em weak abelian fibration} consists
of a base variety $S$, a projective map $f:M\to S$, and a smooth abelian group-scheme $g:P\to S$,
with connected fibers,
acting on $M$.
Thus for $s\in S$, the fiber $M_s = f^{-1}(s)$ is a projective variety, the fiber $P_s = f^{-1}(s)$ is
a connected abelian group, and we have a $P_s$-action on $M_s$.

We require the following properties to hold:
\begin{enumerate}
\item For each $s\in S$, the fibers $M_s$ and $P_s$ have the same dimension.
\item For each $s\in S$, and $m\in M_s$, the stabilizer $Stab_{P_s}(m)\subset P_s$ is affine. 
\item  $P$ has a polarizable Tate module.
\end{enumerate}
\end{defn}

\begin{remark}
We will not attempt to explain the third technical condition other than the following brief remark.
For each $s\in S$, there is a canonical Chevalley exact sequence 
$$
\xymatrix{
1 \ar[r] & R_s \ar[r] & P_s \ar[r] & A_s \ar[r] & 1
}
$$
where $A_s$ is an abelian variety,  and $R_s$ is a connected abelian affine group.
If $R_s$ were always trivial, the third condition would assert that $P$ is a polarizable family of abelian varieties.
\end{remark}

As the name suggests, the above notion is quite weak. To strengthen it, let us cut up the base variety $S$ into its
subvarieties $S_\delta$ where the 
dimension $\delta(s) = \dim(R_s)$ of the affine part of $P_s$ is precisely $\delta$. 

\begin{defn}  
(1) A smooth connected abelian group scheme $g : P \to S$ is said to be
{\em $\delta$-regular} if it satisfies
$$
\on{codim}(S_\delta) \geq  \delta.
$$

(2) A {\em $\delta$-regular abelian fibration} is a weak abelian fibration whose group scheme is $\delta$-regular.
\end{defn}

\begin{remark}
Given any $\delta$-regular abelian fibration, over the generic locus $S_0\subset S$, the group-schemes $P_s$ are in fact abelian varieties, act on $M_s$
with finite stabilizers, and hence $M_s$ is a finite union of abelian varieties.
\end{remark}

\begin{example}
Here are two good examples of $\delta$-regular abelian fibrations to keep in mind:

(1) Any integrable system (though here we use that we are working over the complex numbers $\BC$).

(2) For $X\to S$ a versal deformation of a curve with plane singularities, $M = \ol{Jac}(X/S)$ with its natural action
of $P = Jac(X/S)$.
\end{example}

Now we arrive at the main new technical result underlying Ng\^o's proof of the Fundamental Lemma. It is a refinement
of the celebrated Decomposition Theorem of Beilinson-Bernstein-Deligne-Gabber in the setting of abelian fibrations. In citing the
Decomposition Theorem, we are
invoking the full power of Hodge theory, in particular a very general form of the relative Hard Lefschetz Theorem.
Let us recall it in a specific form sufficient  for our current discussion.

\begin{thm}[Decomposition Theorem]
Let $f:M\to S$ be a projective map of varieties with $M$ smooth. The pushforward $Rf_*\BC_M$ is a direct sum
of (shifted) intersection cohomology sheaves of local systems on subvarieties of $S$.
\end{thm}

\begin{example}

(1) Let $f:M\to S$ be  a proper  fibration with $M$ and $S$ smooth. Then the pushforward $Rf_*\BC_M$
is a direct sum of (shifted) semisimple local systems 
whose fiber at $s\in S$ is the cohomology of the fiber $M_s = f^{-1}(s)$.

(2) Recall the Springer resolution $\mu_{\CN}:\tilde\CN \to \CN$ discussed above. The pushforward $R\mu_{\CN*}\BC_{\tilde\CN}$ is a direct sum of intersection cohomology sheaves supported on the various nilpotent orbits. For simplicity, let us restrict to $G=GL(n)$. Then the nilpotent orbits $\CO_y \subset \CN$ are indexed by Young diagrams $y$. Irreducible representations $V_y$ of the symmetric group $\Sigma_n$
are also indexed by Young diagrams $y$. On each orbit $\CO_y\subset \CN$,  there is a local system $\CL_y$ of rank $\dim V_y$
such that its intersection cohomology is the contribution to  $R\mu_{\CN*}\BC_{\tilde\CN}$  from  the orbit $\CO_y$.
This gives a geometric decomposition of the regular representation of $\Sigma_n$.

\end{example}

The following support result significantly constrains the supports of  the summands occurring in the Decomposition Theorem. Very roughly speaking, it says
that in the case of $\delta$-regular abelian fibrations, one can surprisingly  calculate the pushforward by studying generic loci.

\begin{thm}[Ng\^o's Support Theorem]
 Let $f : M \to S$, $ g : P \to S$ be a $\delta$-regular abelian fibration of relative dimension $d$. Assume the base $S$ is connected and the total space $M$ is smooth.
 
 Let $\CF$ be an intersection cohomology sheaf occurring as a summand in the pushforward  $Rf_*\BC_M$,
and let $Z\subset S$ be the support of $\CF$. 
 
Then there exists an open subset $U \subset S$ such that $U\cap Z\not = \emptyset$,
 and a non-zero local system $\CL$ on the intersection $U \cap Z$ such that the tautological extension by zero of $\CL$
to all of $U$
 is a summand of the restriction $R^{2d} f_* \BC_M|_U$.
 
 In particular, if the fibers of $f$ are irreducible,  then $Z = S$.
\end{thm}

\begin{example} Here is an example and then a non-example of the kind of phenomenon explained by the support theorem.

(1) Let $f:M \to S$ be a proper flat family of irreducible curves with $M$ and $S$ smooth.
Let $S^{reg}\subset S$ be the open locus over which $f$ is smooth, so in particular where the fibers are necessarily smooth.
Then the pushforward $Rf_*\BC_M$ is the intersection cohomology extension of the (shifted) local system $Rf_*\BC_M|_{S^{reg}}$
whose fiber at $s\in S^{reg}$ is the cohomology $H^0\oplus H^1[-1] \oplus H^2 [-2]$ of the curve $M_s = f^{-1}(s)$. 

This satisfies the conclusion of the support theorem, though strictly speaking it is not an application of it 
(except in special instances where the curves 
are of genus $1$).

(2) Consider the proper flat family of irreducible surfaces 
$$
\xymatrix{
f:M = \{ ([x, y, z, w], s) \in \BP^3 \times \BA^1 | x^3 + y^3 + z^3 + sw^3 = 0\} \ar[r] & 
S = \BA^1 
}
$$
with the obvious projection
$
f([x, y, z, w], s) = s.
$
It is smooth over $\BA^1 \setminus\{0\}$, but singular over $\{0\}$. One can check that the pushforward 
$Rf_*\BC_M$ contains summands which are skyscraper sheaves supported at $\{0\}$.

This does not satisfy the conclusion of the support theorem, though  $f$  is a proper flat map.

\end{example}


\subsection{Geometric elliptic endoscopy}
Now let us return to the Hitchin fibration 
along with its relative symmetry group 
$$
\xymatrix{
\chi:\CM_G \ar[r] & \CA_G
&
\CP \ar[r] & \CA_G
}
$$

The results to follow are strikingly parallel to the stabilization of the anisotropic part of the Arthur-Selberg trace formula.
Recall that the ad\`elic factorization of Formula~\eqref{adelic factorization} expressed the anisotropic terms of the Arthur-Selberg trace formula in terms of local orbital integrals, and ultimately led to the Fundamental Lemma. Here we
reverse the process and thanks to the ad\`elic factorization of Theorem~\ref{geometric adelic factorization}, calculate the cohomology of affine Springer fibers in terms of anisotropic Hitchin fibers, and in this way ultimately prove the Fundamental Lemma.

\begin{defn}[Anisotropic Hitchin fibration]
The {anisotropic} Hitchin base 
$\CA_G^{ani} \subset \CA_G$ consists of parameters $a\in \CA_G$ such that any Higgs bundle $(\CP, \varphi)\in \chi^{-1}(a)$ is generically anisotropic in the sense that for any generic trivialization of $\CP$ (and the line bundle $\CL$),
the Higgs field $\varphi$ is anisotropic as an element of $\fg(F)$, where $F$ is the function field of $C$.

The anisotropic Hitchin fibration and relative symmetry group is the restriction of the Hitchin fibration and its relative symmetry group
to the anisotropic Hitchin base
$$
\xymatrix{
\chi^{ani}:\CM_G^{ani} =  \CM_G \times_{\CA_G} \CA_G^{ani} \ar[r] & \CA_G^{ani}
&
\CP^{ani} = \CP\times_{\CA_G} \CA^{ani}_G \ar[r] & \CA^{ani}_G
}
$$
\end{defn}

\begin{remark}
An equivalent characterization of the anisotropic Hitchin base $\CA_G^{ani}$ is the complement of the image of the Hitchin bases for all Levi subgroups.
\end{remark}

\begin{remark}
Since $GL(n)$ has nontrivial split center, it contains no anisotropic tori. Thus the anisotropic Hitchin fibration for $GL(n)$ is empty
and all that follows is vacuous. This is not surprising: there is no nontrivial elliptic endoscopy for $GL(n)$. All of its endoscopic subgroups are in fact Levi subgroups. 

The reader interested in examples is recommended to look at the end of this section where $SL(2)$ is discussed.
\end{remark}

\begin{thm}
Over the anisotropic locus, the Hitchin fibration and its relative symmetry group
%
form a $\delta$-regular abelian fibration.

\end{thm}

\begin{remark}
(1)
In fact, in general $\CM_G^{ani}$ will be a Deligne-Mumford stack, but this does not obstruct any aspect of our discussion.

(2) The theorem is not known in characteristic $p$, but one can restrict further to where it holds and continue
with the calculations  to be performed. Then local-global compatibility can be used to extend the calculations further.
\end{remark}

To shorten the notation, we will write 
$$
\CF_G = R\chi^{ani}_*\BC_{\CM^{ani}_G}
$$
for the pushforward of the constant sheaf along the Hitchin fibration.

Since cohomology is invariant under isotopy,
the natural $\CP^{ani}$-action  on $\CF_G$
factors through the component group 
$\pi_0(\CP^{ani})$.
Our aim consists of two steps:

\begin{enumerate}
\item {\em Fourier theory}: decompose the pushforward 
$\CF_G
$ into eigenspaces
for the natural  $\pi_0(\CP^{ani})$-action,
\item {\em Endoscopic stabilization}: identify the eigenspaces for nontrivial characters with the invariant eigenspaces  for endoscopic groups.
\end{enumerate}
Because we are over the anisotropic locus, $\pi_0(\CP^{ani})$ has finite fibers, and
so the eigenspace decomposition is discrete.

For the trivial character of  $\pi_0(\CP^{ani})$, Ng\^o's Support Theorem gives a striking analogue of the central Theorem~\ref{springer theory}  of Springer theory.

\begin{thm}[Stable summand]
The $\pi_0(\CP^{ani})$-invariant (shifted) simple perverse summands of the pushforward 
$
\CF_G
$ are the canonical intersection cohomology
extensions to all of $\CA_G^{ani}$ of their restrictions to any non-empty open subset.
\end{thm}

Without loss of applicability, we can find an \'etale base change $\tilde \CA^{ani}_G\to \CA^{ani}_G$ over which 
there is a surjective homomorphism of relative groups
$$
\xymatrix{
X_*(T) \ar[r] & \pi_0(\CP^{ani}).
}
$$
This is convenient in that we can think of a character 
$$
\xymatrix{
\kappa:\pi_0(\CP^{ani}) \ar[r] &  \BC^\times
}$$ 
as in turn a character of $X_*(T)$, or in other words, an element of $T^\vee$.

On the one hand, we can consider the summand 
$$
\CF_{G, \kappa} \subset \CF_G
$$ 
consisting of  sections that transform by
the character $\kappa$
under the action of $\pi_0(\CP^{ani})$.
In particular, when $\kappa = e$ is the trivial character, we denote
 the stable summand of $\pi_0(\CP)$-invariants by
$$
\CF_{G, st} = \CF_{G, e}.
$$ 

On the other hand, we can consider endoscopic groups $H$ associated to endoscopic data $(\kappa, \rho_\kappa)$
where $\kappa\in T^\vee\subset G^\vee$ and $\rho_\kappa:\pi_1(C, c_0) \to \pi_0(\kappa)$. Though we are now in a global geometric setting, the constructions are completely parallel to those discussed in Section~\ref{sect: endo}. The resulting endoscopic groups 
$H$ are unramified
quasi-split group schemes defined over $C$. It is straightforward to generalize the Hitchin total space, base, and fibration
to such group schemes. In particular,  the natural transfer map induces a closed immersion of Hitchin bases
$$
\xymatrix{
\nu_H:\tilde\CA^{ani}_H \ar@{^(->}[r] & \tilde\CA^{ani}_G.
}$$

\begin{thm}[Endoscopic summands]
The geometric Fundamental Lemma holds: the
$\kappa$-eigenspace
$
\CF_{G, \kappa}
$
is a direct sum over all $\kappa$-endoscopic groups $H$ of their (shifted) stable summands
$
\nu_{H*}\CF_{H, st}.
$
\end{thm}

The theorem is proved by an intricate application of Ng\^o's support theorem. Roughly speaking,
it is straightforward to establish an identification over a generic open locus intersecting the Hitchin bases of all relevant endoscopic groups. Then one uses the support theorem to conclude that the identification extends over the entire anisotropic Hitchin base.

The identification of  sheaves gives an identification of their fibers  and hence an identification of the stable cohomology of Hitchin fibers for $H$ with the corresponding endoscopic cohomology of Hitchin fibers for $G$.
As we have discussed, such identifications imply analogous identifications for affine Springer fibers, hence also for $p$-adic orbital integrals,
and ultimately confirm the Fundamental Lemma.

 \begin{example}[Geometric endoscopy for $SL(2)$] Recall that to a parameter $a\in \CA_{SL(2)}$,
 we can assign the spectral curve $C_a$
  given by the equation $t^2+a$ on the total space of $\CL$.
As long as $a$ is not identically  zero,  it is generically regular and $C_a$ is reduced.
In this case, we write $\tilde C_a$ for the normalization of $C_a$.


We can cut up $\CA_{SL(2)} \setminus \{0\}$ into three natural pieces:

\begin{enumerate}
\item  $\CA^{st}$: $C_a$ is irreducible, $\tilde C_a\to C$ is ramified. Then $\pi_0(\CP_a)$ is trivial.
\item  $\CA^{U(1)}$: $C_a$ is irreducible, $\tilde C_a\to C$ is unramified. Then $\pi_0(\CP_a) \simeq \BZ/2\BZ$.
\item  $\CA^{GL(1)}$: $C_a$ is reducible. Then $\pi_0(\CP_a) \simeq \BZ$.
\end{enumerate}

The anisotropic locus $\CA_{SL(2)}^{ani}$ is the union $\CA^{st} \sqcup  \CA^{U(1)}$.

Each component of $\CA^{U(1)}$, $\CA^{GL(1)}$ is the image of the Hitchin base for an endoscopic $U(1)$, $GL(1)$ 
respectively.

 \end{example}

\section{Some further directions}\label{sectfuture}

In this section, we briefly list some research directions related to Ng\^o's proof of the Fundamental Lemma.
In particular, we focus on geometric questions, some solved, some open.
It goes without saying that the list is idiosyncratic and far from comprehensive.

\subsection{Deeper singularities}
From a geometric perspective, Ng\^o's endoscopic description of the cohomology of the anisotropic Hitchin fibration is only a first step. In the spirit of Springer theory, one should study the entire fibration, proceeding from the anisotropic locus to  more complicated Higgs bundles. 

Chadouard and Laumon \cite{chla, chlaI, chlaII} have extended Ng\^o's picture to the locus of regular semisimple Higgs bundles in order to prove Arthur's weighted fundamental lemma. Recall that the anisotropic locus consists of Higgs bundles such that the generic value of the Higgs field
is an anisotropic element of $\fg(F)$, where $F$ is the function field of the curve  $C$.
A Higgs bundle is generically regular semisimple if the 
 generic value of the Higgs field
is a regular semisimple element of $\fg(F)$. The regular semisimple Hitchin fibers are no longer of finite type, but are 
 increasing unions of finite-type schemes of bounded dimension. (One could compare with affine Springer fibers which display analogous behavior.) Chadouard and Laumon develop a beautiful truncation framework, directly inspired by Arthur's truncations of the trace formula, and exhibit the regular semisimple Hitchin fibration as an increasing union of proper maps.
 The truncations are governed by stability conditions on Higgs bundles, though the eventual contributions to the weighted fundamental lemma are independent of the stability parameters.

\subsection{Transfer for relative trace formulas}
The relative trace formula\footnote{I would like to thank Y. Sakellaridis for explaining
 to me what this world is all about.} (as described in~\cite{lapid}) is a yet to be fully developed framework,
originally  pioneered by Jacquet, for understanding period integrals of  automorphic forms. It includes a transfer principle
 depending upon identities analogous to the Fundamental Lemma.
 
From a geometric perspective, the main object of study in the relative trace formula is a reductive group $G$ equipped with two spherical subgroups $H_1, H_2\subset G$, and ultimately the two-sided quotient $H_1\bs G/H_2$. The basic example is $G = H \times H$ with $H_1 = H_2 = H$ which gives  the adjoint quotient $G/G$ and the usual trace formula.  

Given two such groups with spherical subgroups $H_1 \times H_2 \subset G, H'_1 \times H'_2 \subset G'$, the transfer principle is predicated on a relation between the $H_1\times H_2$-invariant functions on $G$
 and $H'_1\times H'_2$-invariant functions on $G'$. This  induces a transfer of invariant distributions
   generalizing that between groups and their endoscopic groups
mediated by the relation of their characteristic polynomials.

Geometric methods have been used by Ng\^o~\cite{ngorelative} and Yun~\cite{yunrelative}  to prove fundamental lemmas for relative trace formulas. The latter employs direct analogues of Ng\^o's global proof of the Fundamental Lemma.
There are also precursors to further such results coming from 
traditional Springer theory as in the work of Grinberg~\cite{grinberg}.

\subsection{Geometric trace formulas}
A striking aspect of Ng\^o's proof of the Fundamental Lemma is that its final global calculations are direct geometric analogues of the stabilization of the trace formula. Recall that at the end of the day, the Fundamental Lemma is simply a tool in the analysis of the trace formula and its stabilization.  With the Geometric Langlands program in mind, this naturally leads to the  question: what is the geometric analogue of the trace formula itself? 

The Geometric Langlands program, pioneered by Beilinson and Drinfeld, is a fruitful geometric analogue of the theory of automorphic forms and Galois representations. The basic automorphic objects are
$\mathcal D$-modules on the moduli of principal $G$-bundles on a curve. The corresponding spectral objects are quasicoherent sheaves on
the moduli of flat $G^\vee$-connections on the curve.  
Langlands's conjectural reciprocity takes the form of a conjectural equivalence of categories between $\mathcal D$-modules and quasicoherent sheaves. Broadly understood, the subject forms a substantial industry with motivations from mathematical physics and representation theory.

 Frenkel, Langlands, and Ng\^o~\cite{frenkelngolanglands, frenkelngo} have taken the first steps towards  a trace formula in this setting. They describe the intricate contours of the problem, and make serious progress towards the development of a precise formulation.
  Stepping back from the challenges of principal bundles on curves, one can ask what kind of math should a geometric trace formula involve? Or slightly more precisely, what kind of object is the character of a group acting on a category? Here algebraic topology provide
  precise answers involving loop spaces, and  Hochschild and cyclic homology. For a realization of these ideas in the context of group actions, and in particular, a connection to Lusztig's character sheaves,
   see~\cite{bzn}.

\subsection{Purity of Hitchin fibration}

The decisive geometric input to Ng\^o's proof of the Fundamental Lemma is the Decomposition Theorem applied to the Hitchin fibration. The purity of the pushforward along the Hitchin fibration leads to the endoscopic decomposition
of the cohomology.

If one restricts to equivalence classes of appropriate  semistable Higgs bundles, the Hitchin space takes the more concrete form of a
quasiprojective variety $\CM_{Dol}$. 
Nonabelian Hodge theory provides a diffeomorphism between $\CM_{Dol}$ and a corresponding affine variety $\CM_B$  of  representations of the fundamental group of the curve.
The diffeomorphism is far from an isomorphism: the algebraic structure on $\CM_{Dol}$ depends on the algebraic structure of the curve, while $\CM_B$ depends only on the fundamental group of the curve. Though $\CM_B$ is affine, the fibers of the Hitchin fibration for $\CM_{Dol}$ 
are compact half-dimensional subvarieties.

De Cataldo, Hausel, and Migliorini~\cite{dCHM} study the weight and perverse filtrations on the cohomology of $\CM_{Dol}$ and $\CM_B$ induced by  the Hitchin fibration.
Their specific results and what they point towards should shed further light on the topological nature of the 
indispensable purity invoked in  the proof of the Fundamental Lemma.

\subsection{Affine Springer theory}

The broad paradigms of Springer theory explain many aspects of Ng\^o's proof of the Fundamental Lemma.
For example, it employs in an essential way the idea that the cohomology of complicated Springer fibers could be recovered from  simpler fibers.

But conversely, there are many other aspects of Springer theory which would be worth pursuing in the setting of the Hitchin fibration. For example, Yun~\cite{yunI, yunII, yunIII} studies a tamely ramified version of the Hitchin fibration consisting of Higgs bundles equipped with a compatible flag at a point of the curve. This leads to a generalization of the commutative symmetries appearing in the endoscopic decomposition of the cohomology. Namely, the cohomology becomes a module over the double affine Hecke algebra, and other intriguing relations with quantum algebra appear. 

It would be interesting to find other important aspects of traditional Springer theory: the role of the nilpotent cone,  a resolution of nilpotent Higgs bundles, and the role of the Fourier transform, to name a few.

%
%
%
%
%
%
%
%
%
%
%
%
%
%



%


\end{document}